\documentclass{elsart}

\usepackage{amssymb}

\newcommand{\R}{{\mathbb R}}
\newcommand{\N}{{\mathbb N}}
\newcommand{\Z}{{\mathbb Z}}
\newcommand{\tN}{\widetilde{N}}
\newcommand{\tu}{\widetilde{u}}
\newcommand{\ut}{\underline{\tau}}
\newcommand{\sch}{{\mathcal S}(\R^d)}
\newcommand{\Sp}{\mbox{\textnormal{Sp }}}
\newcommand{\Op}{\mbox{\textnormal{Op}}}
\newcommand{\cst}{\mbox{\textnormal{Cst }}}
\newcommand{\modj}{\langle\xi\rangle}
\newcommand{\fourinv}{\check}
\newcommand{\ind}{{2[\frac{d}{2}]+2}}
\newcommand{\Ind}{n+2+[\frac{d}{2}]+d}
\newcommand{\slf}{\sigma_{lf}}
\newcommand{\sI}{\sigma_I}
\newcommand{\stI}{\sigma_{{\widetilde{I}}}}
\newcommand{\sII}{\sigma_{II}}
\newcommand{\sR}{\sigma_R}
\newcommand{\dsp}{\displaystyle}

\newtheorem{exam}{Example}
\newtheorem{nota}{Notation}
\newtheorem{theo}{Theorem}

\newenvironment{proof}{\noindent \textbf{Proof.} \newline}{
\begin{flushright} $\blacksquare$ \end{flushright} }

\newcommand{\reseteq}{\setcounter{equation}{0}}

\begin{document}

\begin{frontmatter}



\title{Sharp estimates for pseudodifferential operators 
	with symbols of limited 
	smoothness 
	and commutators}


\author{David Lannes}

\address{MAB, Universit\'e Bordeaux I et CNRS, UMR 5466\\
	351 Cours de la Lib\'eration,\\
	33405 Talence Cedex, France.}
\ead{lannes@math.u-bordeaux1.fr}

\begin{abstract}
We consider here pseudo-differential operators whose symbol $\sigma(x,\xi)$
is not infinitely smooth with respect to $x$. Decomposing such symbols into
four -sometimes five- components and using tools of paradifferential
calculus, we derive sharp estimates on the action of such pseudo-differential
operators on Sobolev spaces and give explicit 
expressions for their operator norm
in terms of the symbol $\sigma(x,\xi)$. We also study commutator estimates
involving such operators, and generalize or improve the so-called Kato-Ponce
and Calderon-Coifman-Meyer estimates in various ways.
\end{abstract}

\begin{keyword}
pseudodifferential operators \sep paradifferential calculus \sep 
commutator estimates
\PACS 
\end{keyword}
\end{frontmatter}

%%%%%%%%%%%%%%%%%%%%%%%%%%%%%%%%%%%%%%%%%%%%%%%%%%%%%%%%%%%%%%%%%%%%%%%
%%%%%%%%%%%%%%%%%%%%%%%%%%%%%%%%%%%%%%%%%%%%%%%%%%%%%%%%%%%%%%%%%%%%%%%
%%%%%%%%%%%%%%%%%%%%%%%%%%%%%%%%%%%%%%%%%%%%%%%%%%%%%%%%%%%%%%%%%%%%%%%
\reseteq
\section{Introduction}
%%%%%%%%%%%%%%%%%%%%%%%%%%%%%%%%%%%%%%%%%%%%%%%%%%%%%%%%%%%%%%%%%%%%%%%
%%%%%%%%%%%%%%%%%%%%%%%%%%%%%%%%%%%%%%%%%%%%%%%%%%%%%%%%%%%%%%%%%%%%%%%
%%%%%%%%%%%%%%%%%%%%%%%%%%%%%%%%%%%%%%%%%%%%%%%%%%%%%%%%%%%%%%%%%%%%%%%

%%%%%%%%%%%%%%%%%%%%%%%%%%%%%%%%%%%%%%%%%%%%%%%%%%%%%%%%%%%%%%%%%%%%%%%
%%%%%%%%%%%%%%%%%%%%%%%%%%%%%%%%%%%%%%%%%%%%%%%%%%%%%%%%%%%%%%%%%%%%%%%
\subsection{General setting and description of the results}
%%%%%%%%%%%%%%%%%%%%%%%%%%%%%%%%%%%%%%%%%%%%%%%%%%%%%%%%%%%%%%%%%%%%%%%
%%%%%%%%%%%%%%%%%%%%%%%%%%%%%%%%%%%%%%%%%%%%%%%%%%%%%%%%%%%%%%%%%%%%%%%

Among the
widely known properties of pseudo-differential operators with symbol
in Hormander's class $S^m_{1,0}$, two are discussed in this paper. The first
one concerns their action on Sobolev spaces, and the second one deals with
the properties of commutators.

It is a classical result that for all $\sigma^1\in S^{m_1}_{1,0}$, the operator $\Op(\sigma^1)$ maps $H^{s+m_1}(\R^d)$ into $H^s(\R^d)$
for all $s\in\R$. Moreover, the proof shows that
\begin{equation}
	\label{intro1}
	\forall u\in H^{s+m_1}(\R^d),\qquad
	\vert \Op(\sigma^1)u\vert_{H^s}\leq \underline{C}(\sigma^1)
	\vert u\vert_{H^{s+m_1}}.
\end{equation}

Concerning the study of commutators, Taylor (following works of Moser \cite{Moser} and Kato-Ponce \cite{Kato-Ponce}) proved in \cite{TaylorM0} that for
all $\sigma^1\in S^{m_1}_{1,0}$, with $m_1>0$, and all $\sigma^2\in H^\infty(\R^d)$, one has
for all $s\geq 0$,
\begin{equation}
	\label{intro2}
	\big\vert [\Op(\sigma^1),\sigma^2]u\big\vert_{H^s}
	\leq \underline{C}(\sigma^1)\big(\vert\sigma^2\vert_{W^{1,\infty}} 
	\vert u\vert_{H^{s+m_1-1}}
	+\vert \sigma^2\vert_{H^{s+m_1}}\vert u\vert_\infty\big);
\end{equation}
another well-known commutator estimate is the so-called Calderon-Coifman-Meyer
estimate: if $m_1\geq 0$ then 
for all $s\geq 0$ and $t_0>d/2$ such that $s+m_1\leq t_0+1$, one has
(see Prop. 4.2 of \cite{TaylorM} for instance):
\begin{equation}
	\label{intro3}
	\big\vert [\Op(\sigma^1),\sigma^2]u\big\vert_{H^s}
	\leq \underline{C}(\sigma^1)\vert \sigma^2\vert_{H^{t_0+1}}
	\vert u\vert_{H^{s+m_1-1}}.
\end{equation}
A drawback of (\ref{intro1}), (\ref{intro2}) and (\ref{intro3})
is that the dependence
of the constant $\underline{C}(\sigma^1)$ on $\sigma^1$ is not specified.
This may cause these estimates to be inoperative in the study of some 
nonlinear PDE; indeed, when solving such an equation by an iterative scheme,
one is led to study the pseudo-differential operator corresponding to
the linearized equations around some reference state. Generally, the
symbol of this operator can be written $\sigma(x,\xi)=\Sigma(v(x),\xi)$,
where $\Sigma(v,\xi)$ is smooth with respect to $v$ and of order $m$ with
respect to $\xi$, while $v(\cdot)$ belongs to some Sobolev space $H^s(\R^d)$.
For instance, in the study of nonlinear water waves, one is led to study
the operator associated to the symbol (see \cite{Lannes})
\begin{equation}
	\label{symbDN}	
	\sigma(x,\xi):=\sqrt{(1+\vert\nabla a\vert^2)\vert\xi\vert^2-
	(\nabla a\cdot\xi)^2},
\end{equation}
which is of the form described above, with 
$\Sigma(v,\xi)=\sqrt{(1+\vert v\vert^2)\vert\xi\vert^2-(v\cdot\xi)^2}$ 
and $v(\cdot)=\nabla a$. Such symbols $\sigma(x,\xi)$ are not 
infinitely smooth with respect to $x$, since their regularity is limited
by the regularity of the function $v$. One must therefore be able to handle
symbols of limited smoothness to deal with 
such situations; moreover, one must be
able to say which norms of $v(\cdot)$ are involved in the constant
$\underline{C}(\sigma^1)$ of (\ref{intro1}), (\ref{intro2}) and 
(\ref{intro3}).\\
But even knowing precisely the way the constants $\underline{C}(\sigma^1)$
depend on $\sigma^1$, the estimates (\ref{intro1}) and (\ref{intro2}) may
not be precise enough in some situations. Indeed, when one has to use, say, a
Nash-Moser iterative scheme, \emph{tame} estimates are needed.
For instance, in such situations, the product estimate $\vert uv\vert_{H^s}
\lesssim \vert u\vert_{H^s}\vert v\vert_{H^s}$ ($s>d/2$) 
is inappropriate and must
be replaced by Moser's tame product estimate $\vert uv\vert_{H^s}\lesssim 
\vert u\vert_\infty\vert v\vert_{H^s}+\vert u\vert_{H^s}\vert v\vert_\infty$
($s\geq 0$).
Obviously, (\ref{intro1}) is not precise enough to contain this latter 
estimate. Part of this paper is therefore devoted to the derivation of
sharper versions of (\ref{intro1}).

In the works dealing with pseudo-differential operators with nonregular
symbols, the focus is generally on the continuity of such
operators on Sobolev or Zygmund spaces 
(see for instance \cite{TaylorM0,TaylorM3}) and not on the derivation of
precise (and tame) estimates. In \cite{Grenier}, E. Grenier gave some 
description of the constants $\underline{C}(\sigma^1)$ in 
(\ref{intro1})-(\ref{intro2})
but his results, though sufficient for his purposes, are far from optimal.
In this article, we aim at proving more precise versions of (\ref{intro1}),
(\ref{intro2}) and (\ref{intro3}), 
and we also give some extensions of these results. Let
us describe roughly some of them:\\
{\bf Action of pseudo-differential operators on Sobolev spaces} 
	(see Corollary \ref{coroII1}). Take a symbol 
$\sigma\in S^m_{1,0}$ of the form $\sigma(x,\xi)=\Sigma(v(x),\xi)$, with
$\Sigma$ as described above and $v\in H^\infty$. Then Moser's tame product
estimate can be generalized to pseudo-differential operators of order $m>0$:
for all $s>0$,
$$
	\vert \Op(\sigma)u\vert_{H^s}
	\lesssim C(\vert v\vert_{\infty})
	(\vert v\vert_{H^{s+m}}\vert u\vert_\infty+\vert u\vert_{H^{s+m}}).
$$
Another estimate which does not assume any restriction on the order $m$ and
also holds for negative values of $s$ is the following: for all 
$t_0>d/2$, one has
\begin{eqnarray*}
	\forall -t_0<s < t_0,&\quad&
	\vert \Op(\sigma)u\vert_{H^s}\lesssim C(\vert v\vert_{\infty})
	\vert v\vert_{H^{t_0}}\vert u\vert_{H^{s+m}},\\
	\forall t_0\leq s, &\quad&
	\vert \Op(\sigma)u\vert_{H^s}
	\lesssim C(\vert v\vert_{\infty})
	(\vert v\vert_{H^{s}}\vert u\vert_{H^{m+t_0}}
	+\vert u\vert_{H^{s+m}}).
\end{eqnarray*}
{\bf Commutator estimates.} In this paper, we give a precise description of the
constant $\underline{C}(\sigma^1)$ which appears in (\ref{intro2}) and 
(\ref{intro3}), and 
generalize these estimates in three directions: 
\begin{itemize}
	\item We control the symbolic
	expansion of the commutator in terms of the Poisson brackets.
	For instance, in the particular case when the symbol 
	$\sigma^1(x,\xi)=\sigma^1(\xi)$ does not depend
	on $x$, we derive the following estimate (see Th. \ref{theoIII3}): 
	if $m_1\in\R$ and $n\in\N$ are such that 
	$m_1>n$, then for all $s\geq 0$, one has
	\begin{eqnarray*}
	\lefteqn{\left\vert [\Op(\sigma^1),\sigma^2]u
	-\Op(\{\sigma^1,\sigma^2\}_n)u
	\right\vert_{H^s}}\\
	&\lesssim&C(\sigma^1)\left(
	\vert\nabla^{n+1}\sigma^2\vert_{{\infty}}
	\vert u\vert_{H^{s+m_1-n-1}}+
	\vert \sigma^2\vert_{H^{s+m_1}}
	\vert u\vert_\infty\right),
	\end{eqnarray*}
	and a precise description of 
	$C(\sigma^1)$ is given; if $\sigma^1(\cdot)$ is regular at the origin,
	we have a more precise version involving only derivatives of 
	$\sigma^2$,
	\begin{eqnarray*}
	\lefteqn{\left\vert [\Op(\sigma^1),\sigma^2]u
	-\Op(\{\sigma^1,\sigma^2\}_n)u
	\right\vert_{H^s}}\\
	&\lesssim&C'(\sigma^1)\left(
	\vert\nabla^{n+1}\sigma^2\vert_{{\infty}}
	\vert u\vert_{H^{s+m_1-n-1}}+
	\vert \nabla^{n+1}\sigma^2\vert_{H^{s+m_1-n-1}}
	\vert u\vert_\infty\right).
	\end{eqnarray*}
	For a similar generalization of (\ref{intro3}), see 
	Th. \ref{theoIII1bis}.
	\item We allow $\sigma^2$ to be a pseudo-differential operator 
	and not only a function (Ths. \ref{theoIII1}, \ref{theoIII2},
	\ref{theoIII1bis},
	\ref{theoIV1} and \ref{theoIV1bis},
	and Corollaries \ref{coroIII3} and \ref{coroIV1});
	\item We give an alternative to (\ref{intro2}) allowing the cases 
	$m_1\leq 0$ and $s<0$ (Ths. \ref{theoIII1} and \ref{theoIV1} and
	Corollaries  \ref{coroIII3} and \ref{coroIV1}); similarly, we
	show that negative values of $s$ and $m_1$ are possible in
	(\ref{intro3}) (see Ths. \ref{theoIII1bis} and \ref{theoIV1bis}).
	For instance, if $\sigma^1$ is a Fourier multiplier of order $m_1\in\R$
	and $\sigma^2$ is of order $m_2\in\R$ with
	$\sigma^2(x,\xi)=\Sigma^2(v(x),\xi)$ and 
	$v\in H^{\infty}(\R^d)^p$ then
	for all $s\in\R$ such that $\max\{-t_0,-t_0-m_1\}<s$ 
	 (with $t_0>d/2$ arbitrary),
	\begin{eqnarray*}
	\lefteqn{\big\vert \big[\sigma^1(D),\sigma^2(x,D)\big]u
	-\Op({\{\sigma^1,\sigma^2\}_n})
	u\big\vert_{H^s}}\\
	& \leq &C(\sigma^1,\vert v\vert_{W^{n+1,\infty}})
	\big(\vert u\vert_{H^{s+m_1+m_2-n-1}}
	+
	\vert v\vert_{H^{(s+m_1\wedge n)_+}}
	\vert u\vert_{H^{m_1+m_2+t_0-m_1\wedge n}}\big).
	\end{eqnarray*}
\end{itemize}

\bigbreak

The above results admit generalization to $L^p$-based Sobolev spaces and Besov
spaces, but we deliberately chose to work with classical $L^2$-based Sobolev 
spaces to ease the readability. We refer the reader interested by this
kind of generalizations to \cite{Yamazaki1,Yamazaki2}, \cite{Marschall} and 
\cite{TaylorM0,TaylorM3} for instance.\\
The methods used to prove the above results rely heavily on Bony's 
paradifferential calculus \cite{Bony} as well as on the works of 
Coifman and Meyer \cite{Meyer,Meyer-Coifman3}.\\
In Section \ref{sectsymb},
we introduce the class of symbols adapted to our study; they consist in
all the symbols $\sigma(x,\xi)$ such that $\sigma(\cdot,\xi)$ belongs
to some Sobolev space for all $\xi$. These symbols are decomposed into four
components, one of them being the well-known paradifferential symbol associated
to $\sigma$, and some basic properties are given.\\
In Section \ref{sectoper}, we study the action on Sobolev (and Zygmund) spaces
of the four components into
which each symbol is decomposed, and give precise estimate on the operator
norm. These results generalize classical results of paraproduct theory
and are in the spirit of \cite{Marschall} and especially \cite{Yamazaki1} 
(but the estimate we give
here are different from the ones given in this latter reference). Gathering
the estimates obtained on each component, we obtain a tame estimate on the
action of the operator associated to the full symbol $\sigma(x,\xi)$.\\
Section \ref{sectcomm} is devoted to the study of commutator estimates. We first give in Prop. \ref{lemmIII1} precise estimates for Meyer's well-known result
on the symbolic calculus for paradifferential operators. In Section \ref{sectfour}, we address the case of commutators between a Fourier multiplier $\sigma^1(D)$ and a pseudo-differential operator $\sigma^2(x,D)$; we study some particular cases, including the case when $\sigma^2(x,\xi)=\sigma^2(x)$ is a function. The case when $\sigma^1(x,D)$ is a pseudo-differential operator (and not only a 
Fourier multiplier) is then addressed in Section \ref{sectfin}.

\bigbreak

Throughout this paper, we use the following notations.\\
\noindent
{\bf Notations.} {\bf i.} For all $a,b\in\R$, we write $a\wedge b:=\max\{a,b\}$;\\
{\bf ii.} For all $a\in \R$, we write $a_+:=\max\{a,0\}$, while
$[a]$ denotes the biggest integer smaller than $a$;\\
{\bf iii.} If $f\in F$ and $g\in G$, $F$ and $G$ being two Banach spaces,
the notation $\vert f\vert_{F}\lesssim \vert g\vert_{G}$ means that
$\vert f\vert_{F}\leq C\vert g\vert_{G}$ for some constant $C$ which does
not depend on $f$ nor $g$.\\
{\bf iv.} Here, $\sch$ denotes the Schwartz space of rapidly decaying functions,
and for any distribution $f\in {\mathcal S}'(\R^d)$, we write respectively
$\widehat{f}$ and $\check{f}$ its Fourier and inverse Fourier transform.\\
{\bf v.} We use the classical notation $f(D)$ to 
denote the Fourier multiplier, namely, 
$\widehat{f(D)u}(\cdot)=f(\cdot)\widehat{u}(\cdot)$.

%%%%%%%%%%%%%%%%%%%%%%%%%%%%%%%%%%%%%%%%%%%%%%%%%%%%%%%%%%%%%%%%%%%%%%%
%%%%%%%%%%%%%%%%%%%%%%%%%%%%%%%%%%%%%%%%%%%%%%%%%%%%%%%%%%%%%%%%%%%%%%%
\subsection{Brief reminder of Littlewood-Paley theory}
%%%%%%%%%%%%%%%%%%%%%%%%%%%%%%%%%%%%%%%%%%%%%%%%%%%%%%%%%%%%%%%%%%%%%%%
%%%%%%%%%%%%%%%%%%%%%%%%%%%%%%%%%%%%%%%%%%%%%%%%%%%%%%%%%%%%%%%%%%%%%%%

We recall in this section basic facts in Littlewood-Paley theory.

\bigbreak

Throughout this article, $\psi\in C_0^\infty(\R^d)$ denotes a smooth bump 
function
such that
\begin{equation}
	\label{LP1}
	\psi(\xi)=1\quad\mbox{ if } \quad \vert\xi\vert\leq1/2 
	\quad\mbox{ and } \quad \psi(\xi)=0
	\quad\mbox{ if } \quad \vert\xi\vert\geq 1,
\end{equation}
and we define $\varphi\in C_0^\infty(\R^d)$ as
\begin{equation}
	\label{LP2}
	\varphi(\xi)=\psi(\xi/2)-\psi(\xi),\qquad\forall\xi\in\R^d,
\end{equation}
so that $\varphi$ is supported in the annulus $1/2\leq \vert\xi\vert\leq 2$,
and one has
\begin{equation}
	\label{LP3}	
	1=\psi(\xi)+\sum_{p\geq 0}\varphi(2^{-p}\xi),\qquad\forall \xi\in\R^d.
\end{equation}
For all $p\in\Z$, we introduce the functions $\varphi_p$, supported in
$2^{p-1}\leq\vert\xi\vert\leq 2^{p+1}$, and  defined as
\begin{equation}
	\label{LP4}
	\varphi_p=0\quad\mbox{if} \quad p<-1,\qquad
	\varphi_{-1}=\psi,\qquad
	\varphi_p(\cdot)=\varphi(2^{-p}\cdot)
	\quad\mbox{if} \quad p\geq 0.
\end{equation}
This allows us to give the classical definition of Zygmund spaces:
\begin{defn}
	Let $r\in\R$. Then $C^r_*(\R^d)$ is the set of all 
	$u\in {\mathcal S}'(\R^d)$
	such that
	$$
	\vert u\vert_{C^r_*}:=
	\sup_{p\geq -1}2^{pr}\vert \varphi_p(D)u\vert_\infty<\infty.
	$$
\end{defn}
\begin{rem}
	We recall the continuous embeddings 
	$H^s(\R^d)\subset C_*^{s-d/2}(\R^d)$, for all $s\in\R$,
	and $L^\infty(\R^d)\subset C^0_*(\R^d)$.
\end{rem}

We now introduce admissible cut-off functions, which play an 
important role in paradifferential theory (\cite{Bony,Meyer} and Appendix B of 
\cite{Metivier-Zumbrun}).
\begin{defn}
	\label{defiLP0}
	A smooth function $\chi(\eta,\xi)$ defined on $\R^d\times\R^d$
	is an admissible cut-off function if and only if:
	\begin{itemize}
	\item There are $\delta_1$ and $\delta_2$ such that 
	$0<\delta_1<\delta_2<1$ and
	\begin{eqnarray*}
	\forall \vert\xi\vert\geq 1/2,\qquad
	\chi(\eta,\xi)&=&1 \quad\mbox{ for }\quad \vert\eta\vert\leq 
	\delta_1\vert \xi\vert,\\
	\forall \vert\xi\vert\geq 1/2,\qquad
	\chi(\eta,\xi)&=&0 \quad\mbox{ for }\quad \vert\eta\vert\geq 
	\delta_2\vert\xi\vert;
	\end{eqnarray*}
	\item For all $\alpha,\beta\in\N^d$, there is a constant 
	$C_{\alpha,\beta}$ such that
	\begin{equation}
	\label{estadm}
	\forall (\eta,\xi)\in \R^{2d},
	\quad\vert\xi\vert\geq 1/2,\qquad
	\left\vert\partial_\xi^\alpha\partial_\eta^\beta 
	\chi(\eta,\xi)\right\vert\leq 
	C_{\alpha,\beta}\modj^{-\vert\alpha\vert-\vert\beta\vert}.
	\end{equation}
	\end{itemize}
\end{defn}
\begin{exam}
	A useful admissible cut-off function is given, for all $N\geq 2$, by
	\begin{equation}
	\label{LP7}
	\chi(\eta,\xi)=\sum_{p\geq -1}\psi(2^{-p+N}\eta)\varphi_p(\xi);
	\end{equation}
	One can check that it satisfies indeed the conditions of Definition
	\ref{defiLP0} with $\delta_1=2^{-N-2}$ and $\delta_2=2^{1-N}$.
\end{exam}

An important property satisfied by admissible cut-off functions is that
$\fourinv\chi(\cdot,\xi)$ and its derivatives with respect to $\xi$ 
enjoy good estimates in $L^1$-norm. The next lemma is a simple 
consequence of the estimates imposed in Definition \ref{defiLP0}; we refer
for instance to Appendix B of \cite{Metivier-Zumbrun} for a proof.
\begin{lem}
	\label{lemmLP0}
	Let $\chi(\eta,\xi)$ be an admissible cut-off function. Then for
	all $\beta\in\N^d$,  
	there exists a constant $C_\beta$ such that
	$$
	\forall \xi\in \R^d,
	\qquad	
	\big\vert \partial_\xi^\beta 
	\fourinv\chi(\cdot,\xi)\big\vert_{L^1(\R^d)}\leq 
	C_\beta\modj^{-\vert\beta\vert}.
	$$
\end{lem}

Finally, we end this section with the classical characterization
of Sobolev spaces (see for instance Th. 2.2.1 of \cite{Chemin} or Lemma 9.4 
of \cite{TaylorM3}).
\begin{lem}
	\label{lemmLP1}
	Let $(u_p)_{p\geq -1}$ be a sequence of ${\mathcal S}'(\R^d)$ such that
	for all $p\geq 0$, $\widehat{u_p}$ is supported in 
	$A2^{p-1}\leq\vert\xi\vert\leq B2^{p+1}$, for some $A,B>0$, and 
	such that $\widehat{u_{-1}}$ is compactly supported.\\
	If, for some $s\in\R$, 
	$\dsp \sum_{p\geq -1}2^{2ps}\vert u_p\vert_2^2<\infty$, then
	$$
	\sum_{p\geq -1}u_p=:u\in H^s(\R^d)\quad\mbox{ and }\quad
	\vert u\vert_{H^s}^2\leq \cst \sum_{p\geq -1}2^{2ps}\vert u_p\vert_2^2.
	$$
	Conversely, if $u\in H^s(\R^d)$ then 
	$$
	\sum_{p\geq -1}2^{2ps}\vert\varphi_p(D)u\vert_2^2\leq\cst 
	\vert u\vert_{H^s}^2.
	$$
\end{lem}

%%%%%%%%%%%%%%%%%%%%%%%%%%%%%%%%%%%%%%%%%%%%%%%%%%%%%%%%%%%%%%%%%%%%%%%
%%%%%%%%%%%%%%%%%%%%%%%%%%%%%%%%%%%%%%%%%%%%%%%%%%%%%%%%%%%%%%%%%%%%%%%
%%%%%%%%%%%%%%%%%%%%%%%%%%%%%%%%%%%%%%%%%%%%%%%%%%%%%%%%%%%%%%%%%%%%%%%
\reseteq
\section{Symbols}\label{sectsymb}
%%%%%%%%%%%%%%%%%%%%%%%%%%%%%%%%%%%%%%%%%%%%%%%%%%%%%%%%%%%%%%%%%%%%%%%
%%%%%%%%%%%%%%%%%%%%%%%%%%%%%%%%%%%%%%%%%%%%%%%%%%%%%%%%%%%%%%%%%%%%%%%
%%%%%%%%%%%%%%%%%%%%%%%%%%%%%%%%%%%%%%%%%%%%%%%%%%%%%%%%%%%%%%%%%%%%%%%

As said in the introduction, we are led to consider nonregular symbols 
$\sigma(x,\xi)$ such that
\begin{equation}
	\label{A0}
	\sigma(x,\xi)={\Sigma}(v(x),\xi),
\end{equation}
where $v\in C^0(\R^d)^p$ for some $p\in \N$, while ${\Sigma}$ is a smooth
function belonging to the class $C^\infty(\R^p,{\mathcal M}^m)$ defined below.
\begin{defn}
	\label{defiA0}
	Let $p\in\N$, $m\in \R$ and let 
	${\Sigma}$ be a function defined over $\R^p_v\times\R^d_\xi$. We say
	that  ${\Sigma}\in C^\infty(\R^p,{\mathcal M}^m)$ if 
	\begin{itemize}
	\item 
	$\Sigma_{\vert_{\R^p\times\{\vert\xi\vert\leq 1\}}}\in 
	C^\infty(\R^p;L^\infty(\{\vert\xi\vert\leq 1\}))$;
	\item For all 
	$\alpha\in \N^p$, $\beta\in \N^d$, there exists a nondecreasing
	function $C_{\alpha,\beta}(\cdot)$ such that
	$$
	\sup_{\xi\in\R^d,\vert\xi\vert\geq 1/4}
	\langle \xi\rangle^{\vert\beta\vert-m}
	\left\vert\partial_v^\alpha\partial_\xi^\beta \Sigma(v,\xi)\right\vert
	\leq C_{\alpha,\beta}(\vert v \vert).
	$$
	\end{itemize}
\end{defn}
\begin{exam}
	One can write the symbol $\sigma(x,\xi)$ given in (\ref{symbDN}) under
	the form
	$\sigma(x,\xi)=\Sigma(\nabla a,\xi)$, with 
	$\Sigma(v,\xi)=
	\sqrt{(1+\vert v\vert^2)\vert \xi\vert^2-(v\cdot \xi)^2}$. 
	One can check that $\Sigma\in C^\infty(\R^d,{\mathcal M}^1)$.
\end{exam}
\begin{rem}
	{\bf i. }We do not assume in this article that the symbols are smooth
	at the origin with respect to $\xi$ (for instance, (\ref{symbDN}) has
	singular derivatives at the origin). This is the reason why 
	the estimate in Definition \ref{defiA0} is
	taken over frequencies away from the origin, 
	namely $\vert\xi\vert\geq 1/4$.\\
	{\bf ii.} When $p=0$, then $C^\infty(\R^p,{\mathcal M}^m)$ coincides
	with the class ${\mathcal M}^m$ of symbols of
	Fourier multipliers of order $m$.
\end{rem}

Let us remark now that if $\sigma(x,\xi)$ is as in (\ref{A0}),
then one can write
$$
	\sigma(x,\xi)=\left[ \sigma(x,\xi)-\Sigma(0,\xi)\right] 
	+\Sigma(0,\xi);
$$
the interest of such a decomposition is that the second term is a simple 
Fourier multiplier while the first one is in $H^s(\R^d)$ if 
$v\in H^s(\R^d)^p$, 
$s>d/2$:
\begin{lem}
	\label{lemmA0}
	Let $p\in\N$, $m\in\R$, $s_0>d/2$ and take $v\in H^{s_0}(\R^d)^p$ and
	$\Sigma\in C^\infty(\R^p,{\mathcal M}^m)$; set also
	$\sigma(x,\xi):=\Sigma(v(x),\xi)$.\\
	Defining $\tau(x,\xi)=\sigma(x,\xi)-\Sigma(0,\xi)$,
	one has $\tau(\cdot,\xi)\in H^{s_0}(\R^d)$ for all $\xi\in\R^d$; 
	moreover:
	\begin{itemize}
	\item  One has 
	 $\tau_{\vert_{\R^d\times\{\vert\xi\vert\leq 1\}}}
	\in L^\infty(\{\vert\xi\vert\leq 1\}; H^{s_0}(\R^d))$;
	\item For all $\beta\in \N^d$ and $0\leq s\leq s_0$,
	$$
	\sup_{\xi\in\R^d,\vert\xi\vert\geq 1/4}
	\langle \xi\rangle^{\vert\beta\vert-m}
	\left\vert\partial_\xi^\beta 
	\tau(\cdot,\xi)\right\vert_{H^s}
	\leq C'_{s,\beta}(\vert v\vert_{L^\infty} \vert)
	\vert v\vert_{H^s},
	$$
	where $C'_{s,\beta}(\cdot)$ is some nondecreasing function depending
	only on the $C_{\alpha,\beta}(\cdot)$, 
	$\vert\alpha\vert\leq [s]+2$, introduced in Definition \ref{defiA0}.
	\end{itemize}
\end{lem}
\begin{proof}
This is a simple consequence of Moser's inequality 
(see e.g. Prop. 3.9 of \cite{TaylorM3}) and the properties of $\Sigma$  set 
forth in Definition \ref{defiA0}.
\end{proof}

The previous lemma motivates the introduction of the following class of 
symbols (see also \cite{Yamazaki1,Marschall} for similar symbol classes).
\begin{defn}
	\label{defiA1}
	Let $m\in\R$ and $s_0>d/2$. A symbol
	$\sigma(x,\xi)$ belongs 
	to the class $\Gamma_{s_0}^m$ if and only if
	\begin{itemize}
	\item One has 
	 $\sigma_{\vert_{\R^d\times\{\vert\xi\vert\leq 1\}}}
	\in L^\infty(\{\vert\xi\vert\leq 1\}; H^{s_0}(\R^d))$;
	\item For all $\beta\in \N^d$, one has
	$$
	\sup_{\xi\in\R^d,\vert\xi\vert\geq 1/4}
	\langle \xi\rangle^{\vert\beta\vert-m}
	\left\vert\partial_\xi^\beta 
	\sigma(\cdot,\xi)\right\vert_{H^{s_0}}<\infty .
	$$
	\end{itemize}
\end{defn}
We now set some terminology concerning the regularity of the symbols at the
origin.
\begin{defn}
	\label{defireg}
	We say that $\Sigma\in C^\infty(\R^p,{\mathcal M}^m)$ is 
	\emph{$k$-regular} at the origin if 	
	$\Sigma_{\vert_{\R^p\times\{\vert\xi\vert\leq 1\}}}\in 
	C^\infty(\R^p;W^{k,\infty}(\{\vert\xi\vert\leq 1\}))$.\\
	Similarly, we say that $\sigma\in \Gamma^m_{s_0}$ is \emph{$k$-regular} 
	at the origin  if $\sigma_{\vert_{\R^d\times\{\vert\xi\vert\leq 1\}}}
	\in W^{k,\infty}(\{\vert\xi\vert\leq 1\}; H^{s_0}(\R^d))$.
\end{defn}
\begin{nota}
	It is quite natural to introduce the seminorms $N_{k,s}^m(\cdot)$ 
	and $M^m_{k,l}(\cdot)$
	defined for all $k,l\in\N$, $s\in\R$ and $m\in \R$ as
	\begin{equation}
	\label{A1}
	N_{k,s}^m(\sigma):=\sup_{\vert\beta\vert\leq k}
	\sup_{\vert\xi\vert\geq 1/4}
	\langle\xi\rangle^{\vert\beta\vert-m}
	\left\vert \partial_\xi^\beta\sigma(\cdot,\xi)\right\vert_{H^s}
	\end{equation}
	and
	\begin{equation}
	\label{A2bis}
	M^m_{k,l}(\sigma):=\sup_{\vert \beta \vert\leq k}
	\sup_{\vert\xi\vert\geq 1/4} \modj^{\vert\beta\vert-m}
	\left\vert\partial_\xi^\beta 
	\sigma(\cdot,\xi)\right\vert_{W^{l,\infty}}.
	\end{equation}
	To get information on the low frequencies, we also define
	\begin{equation}
	\label{A2}
	n_{k,s}(\sigma):=
	\sup_{\vert\beta\vert\leq k,\vert\xi\vert\leq 1}
	\left\vert \partial_\xi^\beta \sigma(\cdot,\xi)\right\vert_{H^s}
	\mbox{ and }
	m_k(\sigma):=\sup_{\vert\beta\vert\leq k,\vert\xi\vert\leq 1}
	\vert\partial_\xi^\beta \sigma(\cdot,\xi)\vert_{\infty}.
	\end{equation}
	Note that $M^m_{k,l}(\sigma)$ and $m_k(\sigma)$ still make
	sense when  $\sigma$ is the symbol of a
	Fourier
	multiplier (i.e. if $\sigma(x,\xi)=\sigma(\xi)$). 
	When $l=0$, we simply
	write $M^m_k(\sigma)$ instead of $M^m_{k,0}(\sigma)$.\\
	Finally, the notation $N_{k,s}^m(\nabla_x^l\sigma)$, $l\in\N$, stands
	for $\sup_{\vert\alpha\vert=l}N^m_{k,s}(\partial_x^\alpha\sigma)$;
	we use the same convention for the other seminorms defined above.
\end{nota}

Associated to the class $\Gamma_s^m$ is the subclass of paradifferential
symbols $\Sigma_s^m$ (in the sense of \cite{Bony,Meyer}, see also 
\cite{Meyer-Coifman3} and Appendix B of \cite{Metivier-Zumbrun}). In the
definition below, the notation $\Sp$ is used to denote the spectrum
of a function, that is, the support of its Fourier transform.
\begin{defn}
	\label{defiA2}
	Let $m\in\R$ and $s_0>d/2$. A symbol
	$\sigma(x,\xi)$ belongs 
	to the class $\Sigma_{s_0}^m$ if and only if
	\begin{itemize}
	\item One has $\sigma\in\Gamma_{s_0}^m$,
	\item There exists $\delta\in(0,1)$ such that
	\begin{equation}
	\label{A3}
	\forall \xi\in\R^d,\qquad
	\Sp \sigma(\cdot,\xi)\subset 
	\{\eta\in\R^d,\vert\eta\vert\leq\delta\vert\xi\vert\}.
	\end{equation}
	\end{itemize}
\end{defn}
\begin{rem}
	If $\psi\in C_0^\infty(\R^d)$ is as in (\ref{LP1}), then the 
	spectral condition (\ref{A3}) implies that for all 
	$\sigma\in\Sigma_s^m$, one has
	\begin{equation}
	\label{A4}
	\forall \xi\in\R^d,\qquad
	\sigma(\cdot,\xi)=
	\sigma(\cdot,\xi)*\left[(2\delta\modj)^d
	\fourinv{\psi}(2\delta\modj\cdot)\right].
	\end{equation}
	It is classical (Bernstein's lemma) 
	to deduce that for all $\alpha,\beta\in\N^d$, one has
	\begin{equation}
	\label{A5}
	\forall \xi\in\R^d,\vert\xi\vert\geq 1/4,\quad
	\left\vert \partial_x^\alpha\partial_\xi^\beta
	\sigma(\cdot,\xi)\right\vert_{L^\infty}
	\leq \cst
	M^m_{\vert\beta\vert}(\sigma)
	\modj^{m-\vert\beta\vert+\vert\alpha\vert},
	\end{equation}
	where  
	$M^m_{\vert\beta\vert}(\cdot)$ is defined in (\ref{A2bis}).
\end{rem}

It is well known that symbols of $\Gamma_s^m$ can be smoothed into 
paradifferential symbols of $\Sigma_s^m$. In order to give a precise 
description of the difference between these two symbols (and of the operator 
associated to it), we split every $\sigma\in \Gamma_s^m$ into four components:
\begin{equation}
	\label{A6}
	\sigma(x,\xi)=\slf(x,\xi)+\sI(x,\xi)+\sII(x,\xi)+\sR(x,\xi),
\end{equation}
with, for some $N\in\N$, $N\geq 4$,
\begin{eqnarray}
	\label{A7}
	\slf(\cdot,\xi)&=&\psi(\xi)\sigma(\cdot,\xi),\\
	\label{A8}
	\sI(\cdot,\xi)&=&\sum_{p\geq -1}\psi(2^{-p+N}D_x)\sigma(\cdot,\xi)
	(1-\psi(\xi))\varphi_p(\xi),\\
	\label{A9}
	\sII(\cdot,\xi)&=&\sum_{p\geq -1}\varphi_p(D_x)\sigma(\cdot,\xi)
	(1-\psi(\xi))\psi(2^{-p+N}\xi),\\
	\label{A10}
	\sigma_R(\cdot,\xi)&=&\sum_{p\geq-1}\sum_{\vert p-q\vert\leq N}
	\varphi_q(D_x)\sigma(\cdot,\xi)	(1-\psi(\xi))\varphi_p(\xi),
\end{eqnarray}
where $\psi\in C_0^\infty(\R^d)$ is a bump function satisfying (\ref{LP1}).\\
We also need sometimes a further decomposition of $\sR$ 
as $\sR=\sigma_{R,1}+\sigma_{R,2}$, with
\begin{equation}
	\label{decompreste}
	\sigma_{R,1}(\cdot,\xi):=(1-\psi(D_x))\sR(\cdot,\xi)
	\quad\mbox{and}\quad
	\sigma_{R,2}(\cdot,\xi):=\psi(D_x)\sR(\cdot,\xi);
\end{equation}
note that $\sigma_{R,2}$ is given by the \emph{finite} sum
$$
	\sigma_{R,2}(\cdot,\xi)=\sum_{p\leq N+1}\sum_{\vert p-q\vert\leq N}
	\psi(D_x)\varphi_q(D_x)\sigma(\cdot,\xi)(1-\psi(\xi))\varphi_p(\xi).
$$
\begin{rem}
	\label{rempara}
	{\bf i.} The fact that the sum of the four terms given in 
	(\ref{A7})-(\ref{A10}) equals $\sigma(x,\xi)$ follows directly from 
	(\ref{LP3}).\\
	{\bf ii.} When $\sigma(x,\xi)=\sigma(x)$ does not depend on $\xi$,
	one has $\Op(\slf)u=\sigma\psi(D)u$, 
	$\Op(\sI)u=T_\sigma\widetilde{u}$, 
	$\Op(\sII)u=T_{\widetilde{u}}\sigma$ and 
	$\Op(\sR)u=R(\sigma,\widetilde{u})$, with $\widetilde{u}:=(1-\psi(D))u$
	and where $T_f$ denotes the usual paraproduct operator associated to
	$f$ and $R(f,g)=fg-T_f g-T_gf$ (see \cite{Bony,Meyer,Chemin}).\\
	{\bf iii.} In \cite{Yamazaki1}, Yamazaki used a similar decomposition
	of symbols into three components. We need a fourth one here, namely
	$\slf$, in order to take into account symbols which are not infinitely
	smooth with respect to $\xi$ at the origin. A fifth component is
	also introduced in (\ref{decompreste}); it is used in the proof of 
	the second parts of Ths. \ref{theoIII1}-\ref{theoIII1bis}.
\end{rem}

In the next proposition, we check that $\sI$ belongs to the class of 
paradifferential symbols $\Sigma_s^m$.
\begin{prop}
	\label{propA1}
	Let $m\in\R$, $s_0>d/2$, and let $\sigma\in\Gamma_{s_0}^m$. \\
	Then, the symbol $\sI$ defined in (\ref{A8}) belongs to 
	$\Sigma_{s_0}^m$
	and, for all $k\in\N$ and $s\leq s_0$,
	$$
	N^m_{k,s}(\sI)\leq \cst N_{k,s}^m(\sigma)
	\quad\mbox{ and }\quad
	M_k^m(\sigma_I)\leq \cst M_k^m(\sigma).
	$$
\end{prop}
\begin{proof}
One can write $\sI(\cdot,\xi)
=(1-\psi(\xi))\fourinv\chi(\cdot,\xi)*\sigma(\cdot,\xi)$, where 
$\chi(\eta,\xi)$ denotes the admissible cut-off function constructed in 
(\ref{LP7}). The spectral property (\ref{A3}) is
thus obviously satisfied by $\sI$ and the result follows therefore from simple 
convolution estimates, 
together with the bounds on the $L^1$-norm on the derivatives 
$\partial_\xi^\alpha \fourinv\chi(\cdot,\xi)$ given in Lemma \ref{lemmLP0}.
\end{proof}

Together with the decomposition given in (\ref{A6}), 
we shall also need another kind of 
decomposition, namely, Coifman and Meyer's decomposition into elementary
symbols. The proof of the next proposition is a quite close adaptation of
the proof of Prop. II.5 of \cite{Coifman-Meyer}; it is given in 
Appendix \ref{appCM}.
\begin{prop}
	\label{propA2}
	Let $m\in\R$ and $s_0>d/2$, and let $\sigma\in \Gamma_{s_0}^m$.
	With $\psi\in C^\infty_0(\R^d)$ as given by (\ref{LP1}), one has
	$$
	(1-\psi(\xi))\sigma(x,\xi)=
	\sum_{k\in\Z^d}\frac{1}{(1+\vert k\vert^2)^{[\frac{d}{2}]+1}}
	p_k(x,\xi)\modj^m,
	$$
	with $\dsp p_k(x,\xi)=\sum_{q\geq -1}c_{k,q}(x)\lambda_k(2^{-q}\xi)$, 
	and where:\\
	{\bf i.} The coefficients $c_{k,q}(\cdot)$ are in $H^{s_0}(\R^d)$
	and for all 
	$s\leq s_0$, one has
	$$
	\vert c_{k,q}\vert_{H^s}\leq \cst N^m_{\ind,s}(\sigma);
	$$
	moreover, for all $p\geq -1$, one has 
	$\dsp\vert \varphi_p(D)c_{k,q}\vert_{H^s}
	\leq\cst N^m_{\ind,s}(\sigma_{(p)})$, where the symbol $\sigma_{(p)}$ 
	is defined as
	$\sigma_{(p)}(\cdot,\xi)=\varphi_p(D_x)\sigma(\cdot,\xi)$.\\
	{\bf ii. } For all $k\in \Z^d$, the functions $\lambda_k(\cdot)$ are
	smooth and supported in $2/5\leq \vert \xi\vert\leq 12/5$. Moreover,
	$\dsp\left\vert\fourinv\lambda_k\right\vert_{L^1}$ is bounded from
	above uniformly in $k\in\Z^d$.
\end{prop}

%%%%%%%%%%%%%%%%%%%%%%%%%%%%%%%%%%%%%%%%%%%%%%%%%%%%%%%%%%%%%%%%%%%%%%%
%%%%%%%%%%%%%%%%%%%%%%%%%%%%%%%%%%%%%%%%%%%%%%%%%%%%%%%%%%%%%%%%%%%%%%%
%%%%%%%%%%%%%%%%%%%%%%%%%%%%%%%%%%%%%%%%%%%%%%%%%%%%%%%%%%%%%%%%%%%%%%%
\reseteq
\section{Operators}\label{sectoper}
%%%%%%%%%%%%%%%%%%%%%%%%%%%%%%%%%%%%%%%%%%%%%%%%%%%%%%%%%%%%%%%%%%%%%%%
%%%%%%%%%%%%%%%%%%%%%%%%%%%%%%%%%%%%%%%%%%%%%%%%%%%%%%%%%%%%%%%%%%%%%%%
%%%%%%%%%%%%%%%%%%%%%%%%%%%%%%%%%%%%%%%%%%%%%%%%%%%%%%%%%%%%%%%%%%%%%%%

To any symbol $\sigma(x,\xi)\in C^0(\R^d\times\R^d\backslash\{0\})$, one can
associate an operator $\sigma(x,D)=\Op(\sigma)$ acting on functions whose
Fourier transform is smooth and compactly supported in $\R^d\backslash\{0\}$:
$$
	\forall x\in \R^d,\qquad
	\Op(\sigma)u(x)=(2\pi)^{-d}\int_{\R^d}
	e^{ix\cdot\xi}\sigma(x,\xi)\widehat{u}(\xi)d\xi.
$$
The aim of this section is to study $\Op(\sigma)$ when $\sigma\in \Gamma^m_s$.
In order to do so,we study successively $\Op(\slf)$, $\Op(\sI)$, $\Op(\sII)$ 
and $\Op(\sR)$, where $\slf$, $\sI$, $\sII$ and $\sR$ are the four components
of the decomposition (\ref{A6}).\\
The operator $\Op(\slf)$ is handled as follows:
\begin{prop}
	\label{propII1}
	Let $m\in \R$ and $s_0>d/2$, and let $\sigma\in \Gamma_{s_0}^m$.\\
	{\bf i.} The operator $\Op(\slf)$ extends as an operator 
	mapping any Sobolev space
	into $H^s(\R^d)$, for all $s\leq s_0$. Moreover,
	$$
	\forall t\in \R,\quad\forall s\leq s_0,\quad \forall u\in H^t(\R^d),
	\qquad
	\left\vert \Op(\slf)u\right\vert_{H^s}\lesssim n_{0,s}(\sigma)
	\vert u\vert_{H^t},
	$$
	where $n_{0,s}(\sigma)$ is defined in (\ref{A2}).\\
	{\bf ii.} If $\sigma$ is $\ind$-regular at the origin, 
	the following estimates also hold,
	for all $s\leq s_0$,
	$$
	\vert \Op(\slf)u\vert_{H^s}
	\lesssim n_{\ind,s}(\sigma)\vert u\vert_{C_*^0}
	\lesssim n_{\ind,s}(\sigma)\vert u\vert_\infty.
	$$
\end{prop}
\begin{proof}
By definition, one has
\begin{eqnarray*}
	\Op(\slf)u(x)&=&(2\pi)^{-d}\int_{\R^d}e^{ix\cdot\xi}
	\psi(\xi)\sigma(x,\xi)
	\widehat{u}(\xi)d\xi\\
	&=&(2\pi)^{-d}\int_{\vert\xi\vert\leq 1}e^{ix\cdot\xi}
	\psi(\xi)\sigma(x,\xi)
	\widehat{u}(\xi)d\xi.
\end{eqnarray*}
Since $\dsp \vert x\mapsto e^{ix\cdot\xi}\sigma(x,\xi)\vert_{H^s}\leq 
	\cst \modj^{\vert s\vert}\vert \sigma(\cdot,\xi)\vert_{H^s}$, one has 
$$
	\vert \Op(\slf)u\vert_{H^s}\leq \cst n_{0,s}(\sigma)
	\int_{\vert\xi\vert\leq 1}\modj^{\vert s\vert}
	\vert \widehat{u}(\xi)\vert d\xi.
$$
One can then obtain the first estimate by a simple 
Cauchy-Schwarz argument.\\
In order to prove the second estimate, remark that a simple expansion in
Fourier series shows that
$$
	\psi(\xi)\sigma(x,\xi)=
	{\mathbf 1}_{\{\vert\xi\vert\leq 1\}}(\xi)\sum_{k\in\Z^d}
	\frac{1}{(1+\vert k\vert^2)^{1+[\frac{d}{2}]}}
	c_k(x)e^{i\xi\cdot k},
$$
where $	{\mathbf 1}_{\{\vert\xi\vert\leq 1\}}$ is the characteristic function
of the ball $\{\vert\xi\vert\leq 1\}$ and
$$
	c_k(x)=(1+\vert k\vert^2)^{1+[d/2]}(2\pi)^{-d}
	\int_{[-\pi,\pi]^d}e^{-i\xi\cdot k}\psi(\xi)\sigma(x,\xi)d\xi.
$$
Using methods similar to those used in the proof of Prop. \ref{propA2}, 
one obtains that  
$\vert c_k(\cdot)\vert_{H^s}\leq n_{\ind,s}(\sigma)$.
Since 
$$
\Op(c_k(x)e^{i\xi\cdot k}{\mathbf 1}_{\{\vert\xi\vert\leq 1\}}(\xi))u
=c_k(\cdot)({\mathbf 1}_{\{\vert\xi\vert\leq 1\}}(D)u)(\cdot+k),
$$ 
the result follows from the next lemma:
\begin{lem}
	\label{lemmII0}
	Let $u,v\in\sch$ and assume that $\widehat{v}$ is supported
	in the ball $\{\vert\xi\vert\leq A\}$, for some $A>0$. 
	Then for all $s\in \R$, one has
	$$
	\vert u v\vert_{H^s}\leq \cst
	\vert u\vert_{H^s}\vert v\vert_\infty.
	$$
\end{lem}
\begin{proof}
Write $uv=\sum_{q\geq -1}v\varphi_q(D)u$; except the first
ones, each term of this sum has its spectrum included in an annulus
of size $\sim 2^q$. Thanks to Lemma \ref{lemmLP1}, the $H^s$-norm of
the product $u v$ can therefore be controlled in terms of
$\vert v\varphi_q(D)u \vert_{L^2}$, $q\geq -1$.
Since these quantities are easily bounded from above by 
$\vert v\vert_\infty\vert \varphi_q(D)u\vert_{L^2}$, the lemma follows 
from another application of
Lemma \ref{lemmLP1}.
\end{proof}
\end{proof}

We now turn to study $\Op(\sI)$. As already said, $\sI$ is the 
paradifferential symbol associated to $\sigma$ so that it is well-known that 
$\Op(\sI)$ maps $H^{s+m}$ into $H^s(\R^d)$ for all $s\in\R$ 
(see \cite{Bony,Meyer} and Prop. B.9 of \cite{Metivier-Zumbrun}). However, 
since we need a precise estimate on the operator norm of $\Op(\sI)$, we 
cannot omit the proof.
\begin{prop}
	\label{propII2}
	Let $m\in\R$ and $s_0>d/2$, and let $\sigma\in \Gamma_{s_0}^m$.\\
	If $\sI$ is as defined in (\ref{A8}), then $\Op(\sI)$ extends as a 
	continuous mapping on $H^{s+m}(\R^d)$ with values in $H^s(\R^d)$,
	for all $s\in\R$. Moreover,
	$$
	\forall s\in\R,\quad\forall u\in H^{s+m}(\R^d),\qquad
	\left\vert \Op(\sI)u\right\vert_{H^s}\lesssim
	M_d^m(\sigma)\vert u\vert_{H^{s+m}},
	$$
	where $M_d^m(\sigma)$ is defined in (\ref{A2bis}).
\end{prop}
\begin{proof}
Let us first prove the following lemma, which deals with the action of operators
whose symbol satisfies the spectral property (\ref{A3}).
\begin{lem}
	\label{lemmII1}
	Let $m\in\R$ and 
	$\sigma(x,\xi)\in  
	C^\infty(\R^d\times \R^d\backslash\{0\})$ be such that 
	$M^m_d(\sigma)<\infty$, where $M^m_d(\cdot)$ is as defined in 
	(\ref{A2bis}).\\
	If moreover $\sigma(x,\xi)$ vanishes for $\vert\xi\vert\leq 1/2$
	and satisfies the spectral condition (\ref{A3})
	then $\Op\left(\sigma\right)$ extends as a continuous 
	mapping on $H^{s+m}(\R^d)$
	with values in $H^s(\R^d)$ for all $s\in\R$ and
	$$
	\forall u\in H^{s+m}(\R^d),\qquad
	\left\vert \Op\left(\sigma\right)u\right\vert_{H^s}\lesssim
	\sup_{\vert\xi\vert\geq 1/2}\sup_{\vert\beta\vert\leq d}
	\Big(\modj^{\vert\beta\vert-m}
	\vert \partial_\xi^\beta\sigma(\cdot,\xi)\vert_\infty\Big) 
	\vert u\vert_{H^{s+m}}.
	$$
\end{lem}
\begin{proof}
Using (\ref{LP3}), we write $\dsp\sigma(x,\xi)=
\sum_{p\geq -1}\sigma_p(x,\xi)$, with 
$\sigma_p(x,\xi)=\varphi_p(\xi)\sigma(x,\xi)$. 
For all $u\in\sch$, (\ref{LP3}) and (\ref{LP4}) yield
\begin{equation}
	\label{B1}
	\Op(\sigma)u=\sum_{p\geq -1}\Op(\sigma_p)
	\sum_{\vert p-q\vert\leq 1}\varphi_q(D)u.
\end{equation}
Let us now define $\widetilde{\sigma_p}(x,\xi):=\sigma_p(2^{-p}x,2^p\xi)$ 
for all $p\in\N$. One obviously has (see e.g. Lemma II.1 of 
\cite{Coifman-Meyer}) $\Vert \Op(\sigma_p)\Vert_{L^2\to L^2}
=\Vert \Op(\widetilde{\sigma_p})\Vert_{L^2\to L^2}$. Moreover, Hwang proved
in \cite{Hwang} that
$$
	\Vert \Op(\widetilde{\sigma_p})\Vert_{L^2\to L^2}\leq 
	\cst \sum_{\alpha,\beta\in\{0,1\}^d}
	\left\vert\partial_x^\alpha\partial_\xi^\beta 
	\widetilde{\sigma_p}\right\vert_{L^\infty(\R^d\times\R^d)};
$$
it follows therefore from (\ref{A5}) that
\begin{equation}
	\label{B2}
	\Vert \Op({\sigma_p})\Vert_{L^2\to L^2}
	\leq \cst M^m_d(\sigma) 2^{pm}.
\end{equation}
The result follows therefore from (\ref{B1}), (\ref{B2}) and 
Lemma \ref{lemmLP1} if we can prove that for all $p\in\N$, 
$\dsp \Op(\sigma_p)\sum_{\vert p-q\vert\leq 1}\varphi_q(D)u$ has its 
spectrum supported in an annulus $A 2^{p-1}\leq \vert\eta\vert\leq B 2^{p+1}$ 
for some $A,B>0$. Since this is an easy consequence
of the spectral property (\ref{A3}), the proof of the lemma is complete.
\end{proof}
The proof of the proposition is now very simple. One just has to apply 
Lemma \ref{lemmII1} to $\sI$, and to use Prop. \ref{propA1}.
\end{proof}

The next proposition gives details on the action of $\Op(\sII)$.
\begin{prop}
	\label{propII3}
	Let $m\in\R$ and $s_0>d/2$, and let $\sigma\in \Gamma_{s_0}^m$.\\
	If $\sII$ is defined as in (\ref{A9}) 
	then $\Op(\sII)$ extends as an operator on any Sobolev space 
	and one has, for all $s\leq s_0$ and $t>0$,
	$$
	\forall u\in C_*^{t+m},\qquad 
	\left\vert\Op(\sII)u\right\vert_{H^s}\lesssim 
	N_{\ind,s}^m(\sigma)\vert u\vert_{C_*^{t+m}},
	$$
	and
	$$
	\forall u\in C_*^{-t+m},\qquad 
	\left\vert\Op(\sII)u\right\vert_{H^{s-t}}\lesssim 
	N_{\ind,s}^m(\sigma)\vert u\vert_{C_*^{-t+m}},
	$$
	where $N^m_{\ind,s}(\sigma)$ is  defined in (\ref{A1}).
\end{prop}
\begin{rem}
	\label{remagrad}
	{\bf i.} One can replace the quantity $N^m_{\ind,s}(\sigma)$ by 
	$N^m_{\ind,s-k}(\nabla_x^k\sigma)$, $k\in\N$, in the estimates
	of the proposition. This follows from the fact that 
	$\vert f\vert_{H^s}\leq \cst \vert \nabla^k f\vert_{H^{s-k}}$, 
	$k\in\N$, whenever $\widehat{f}$ vanishes in a neighborhood of 
	the origin, and from the observation that one can
	replace $\sigma$ by $(1-\psi(D_x))\sigma$ in 
	the definition of $\sII$.\\
	{\bf ii.}
	As said previously, when $\sigma(x,\xi)=\sigma(x)$ does not
	depend on $\xi$, one has $\Op(\sII)u=T_{(1-\psi(D))u}\sigma$ and
	thus $\vert \Op(\sII)u\vert_{H^s}\lesssim 
	\vert u\vert_\infty \vert \sigma\vert_{H^s}$, that is, the endpoint 
	case $t=0$ holds in Prop. \ref{propII3} 
	if one weakens the $\vert u\vert_{C^0_*}$-control
	into a $\vert u\vert_\infty$-control. This is no longer true in general
	when dealing with general symbols.
\end{rem}
\begin{proof}
Prop. \ref{propA2} allows us to reduce the study to the case $m=0$ and
to the reduced symbols $p_k(x,\xi)$ given in that proposition.\\
By definition of $\sII$, one has, for all $u\in\sch$, 
$\sII(x,D)u=\sum_{p\geq-1}v_p$, with 
$v_p=\sigma_{(p)}(x,D)(1-\psi(D))\psi(2^{-p+N}D)u$ and $\sigma_{(p)}(\cdot,\xi)=\varphi_p(D_x)\sigma(\cdot,\xi)$. 
Since the spectrum of $v_p$ is supported in 
$(1-2^{1-N})2^{p-1}\leq \vert\xi\vert\leq (1+2^{-1-N})2^{p+1}$, 
Lemma \ref{lemmLP1} reduces the control of $\vert \sII(x,D)u\vert_{H^s}$ 
to finding an estimate on each $\vert v_p\vert_2$, and hence on
\begin{eqnarray*}
	I&=&\Big\vert \sum_{q\geq -1}\big[\varphi_p(D_x)c_{k,q}\big]
	\lambda_k(2^{-q}D)
	\psi(2^{-p+N}D)(1-\psi(D))u\Big\vert_2\\
	&\leq&\sum_{q\geq -1}\left\vert \varphi_p(D_x)c_{k,q}\right\vert_2
	\left\vert \lambda_k(2^{-q}D)
	\psi(2^{-p+N}D)u\right\vert_\infty.
\end{eqnarray*}
Remarking that $\lambda_k(2^{-q}\xi)\psi(2^{-p+N}\xi)=0$ when $q\geq p-N+2$, 
and using Prop. \ref{propA2}, one deduces
\begin{equation}
	\label{II1}
	I\leq \cst N^m_{\ind,0}(\sigma_{(p)})\sum_{q=-1}^{p-N+1}
	\left\vert \lambda_k(2^{-q}D)
	\psi(2^{-p+N}D)u\right\vert_\infty,
\end{equation}
where we recall that 
$\sigma_{(p)}(\cdot,\xi)=\varphi_p(D_x)\sigma(\cdot,\xi)$.\\
We now need the following lemma:
\begin{lem}
	\label{lemmII2}
	Let $A,B>0$ and $\lambda\in C_0^\infty(\R^d)$ supported in 
	$A\leq\vert\xi\vert\leq B$. Then, for all $t\in\R$ and
	$q\geq -1$, one has, 
	$$
	\forall u\in C_*^t,\qquad
	\left\vert \lambda(2^{-q}D)u\right\vert_\infty
	\leq C_t \vert \fourinv\lambda\vert_{L^1}2^{-qt}\vert u\vert_{C^t_*}.
	$$
\end{lem}
\begin{proof}
Since $\lambda$ is supported in $A\leq\vert\xi\vert\leq B$, there exists 
$n_0\in \N$ such that for all $\xi\in\R^d$ and $q\geq -1$, one has 
$\dsp \lambda(2^{-q}\xi)=\lambda(2^{-q}\xi)\sum_{\vert r-q\vert\leq n_0}
\varphi_r(\xi)$. Therefore, one can write
\begin{eqnarray*}	
	\vert\lambda(2^{-q}D)u\vert_\infty&=&
	\Big\vert \lambda(2^{-q}D)\sum_{\vert r-q\vert\leq n_0}\varphi_r(D)
	u\Big\vert_\infty\\
	&\leq&\vert\fourinv\lambda\vert_{L^1}\sum_{\vert r-q\vert\leq n_0}
	\vert\varphi_r(D)u\vert_\infty,
\end{eqnarray*}
so that the lemma follows from the very definition of Zygmund spaces.
\end{proof}
In order to prove the first part of Prop. \ref{propII3}, take any $t>0$ and
use the lemma to remark that (\ref{II1}) yields
\begin{eqnarray}
	\nonumber
	I&\leq& \cst N^m_{\ind,0}(\sigma_{(p)})\sum_{q\geq -1}
	2^{-qt}\left\vert
	\psi(2^{-p+N}D)u\right\vert_{C^t_*}\\
	\label{bis}
	&\leq& \cst N^m_{\ind,0}(\sigma_{(p)})
	\left\vert u\right\vert_{C_*^{t}},
\end{eqnarray}
since $\sum_{q\geq -1}2^{-qt}<\infty$.
From Lemma \ref{lemmLP1} and the definition of $N^m_{\ind,0}(\cdot)$, it is 
obvious that 
\begin{equation}
	\label{truc}
	\sum_{p\geq -1}2^{2ps}N^m_{\ind,0}(\sigma_{(p)})^2\leq \cst
	N^m_{\ind,s}(\sigma)^2,
\end{equation}
so that (\ref{bis}) and Lemma \ref{lemmLP1} give the result.\\
To prove the second part of the proposition, 
proceed as above to obtain
\begin{eqnarray*}
	I&\leq& \cst N^m_{\ind,0}(\sigma_{(p)})\left(\sum_{q=-1}^{p-N+1}2^{qt}
	\right)\left\vert u\right\vert_{C_*^{-t}},\\
	&\leq&  \cst N^m_{\ind,0}(\sigma_{(p)})2^{pt}
	\left\vert u\right\vert_{C_*^{-t}};
\end{eqnarray*}
the end of the proof is done as for the first part of the proposition.
\end{proof}

We finally turn to study $\Op(\sR)$:
\begin{prop}
	\label{propII4}
	Let $m\in\R$ and $s_0>d/2$, 
	and let $\sigma\in\Gamma^m_{s_0}$.\\
	If $\sR$ is as given in (\ref{A10}) and if
	$s+t>0$ and $s\leq s_0$ then $\Op(\sR)$ extends as
	a continuous operator on $H^{m+t}(\R^d)$ with values in 
	$H^{s+t-\frac{d}{2}}(\R^d)$. Moreover,
	$$
	\forall u\in H^{m+t}(\R^d),\qquad
	\left\vert \Op(\sigma_R)u\right\vert_{H^{s+t-d/2}}
	\lesssim N^m_{\ind,s}(\sigma)\vert u\vert_{H^{m+t}}
	$$
	and
	$$
	\forall u\in C_*^{m+t}(\R^d),\qquad
	\left\vert \Op(\sigma_R)u\right\vert_{H^{s+t}}
	\lesssim N^m_{\ind,s}(\sigma)\vert u\vert_{C_*^{m+t}}.
	$$
\end{prop}
\begin{rem}
	\label{remagrad1}
	For the same reasons as in Remark \ref{remagrad}, one can replace
	the quantity $N^m_{\ind,s}(\sigma)$ by 
	$N^m_{\ind,s-k}(\nabla_x^k\sigma)$, $k\in\N$, in the estimates of the 
	proposition, provided that one replaces
	$\sR$ by $\sigma_{R,1}$, where $\sigma_{R,1}$ is defined in 
	(\ref{decompreste}).
\end{rem}
\begin{proof}
We only prove the first of the two estimates given in the proposition. The 
second one is both easier and contained in 
Th. B of \cite{Yamazaki1}.
The proof we present below is an adaptation of the corresponding result
which gives control of the residual term in paraproduct theory
(e.g. Th. 2.4.1. of \cite{Chemin}).\\
Using the expression of $\sR$ given in (\ref{A10}) and a Littlewood-Paley
decomposition, one can write
$$
	\Op(\sR)u=\sum_{r\geq -1}\varphi_r(D)\Op(\sR)u
	=\sum_{r\geq -1}\varphi_r(D)\sum_{p\geq -1}R_p(\sigma)u,
$$
where $\dsp R_p(\sigma)u:=\sum_{\vert p-q\vert\leq N}
\sigma_{(q)}(x,D)(1-\psi(D))\varphi_p(D)u$, and with
$\sigma_{(q)}(\cdot,\xi)=\varphi_q(D)\sigma(\cdot,\xi)$.\\
Since $\Sp R_p(\sigma)u$ is included in $\vert\xi\vert\leq (1+2^N)2^{p+1}$,
there exists $n_0\in\N$ such that $\varphi_r(D)R_p(\sigma)u=0$ whenever 
$r>p+n_0$. Thus, one has in fact
$$
	\Op(\sR)u=\sum_{r\geq -1}\varphi_r(D)\sum_{p\geq r-n_0}R_p(\sigma)u,
$$
and the proposition follows from Lemma \ref{lemmLP1} and the estimate
\begin{equation}
	\label{II2}
	\Big(\sum_{r\geq -1}2^{2r(s+t-d/2)}
	\Big\vert \varphi_r(D)\sum_{p\geq r-n_0}
	R_p(\sigma)u\Big\vert_2^2\Big)^{1/2}
	\lesssim N^m_{\ind,s}(\sigma)\vert u\vert_{m+t}.
\end{equation}
The end of the proof is thus devoted to establishing (\ref{II2}).\\
Using Prop. \ref{propA2} --and with the same notations-- one can see that
it suffices to prove (\ref{II2}) with $R_p(\sigma)$ replaced by 
$R_p(p_k\modj^m)$, provided that the estimate is uniform in $k\in\Z^d$.
Without loss of generality, we can also assume that $m=0$.\\
Now, remark that
\begin{equation}	
	\label{II3}
	2^{r(s+t-d/2)}\Big\vert \varphi_r(D)\sum_{p\geq r-n_0}
	R_p(p_k)u\Big\vert_2
	\leq \cst 2^{r(s+t)}\Big\vert\sum_{p\geq r-n_0}R_p(p_k)
	u\Big\vert_{L^1}
\end{equation}
and that
\begin{equation}
	\label{II4}
	\left\vert R_p(p_k)u\right\vert_{L^1}=
	\Big\vert \sum_{\vert p-q\vert\leq N}\varphi_q(D_x)p_k(x,D)
	(1-\psi(D))\varphi_p(D) u\Big\vert_{L^1}.
\end{equation}
Now, using the expression of $p_k(x,\xi)$ given in Prop. \ref{propA2}, one can 
write
\begin{eqnarray*}
	\lefteqn{\varphi_q(D_x)p_k(x,D)(1-\psi(D))
	\varphi_p(D)u
	=}\\
	& &\sum_{l\geq -1}\varphi_q(D_x)c_{k,l}(x)\lambda_k(2^{-l}D)
	(1-\psi(D))\varphi_p(D) u,
\end{eqnarray*}
and since $\lambda_k(2^{-l}\xi)\varphi_p(\xi)=0$ if $\vert p-l\vert>n_1$, for
some $n_1\in\N$, one deduces that the summation in the r.h.s. of the above
inequality is over a finite number of integers $l$; therefore, by 
Cauchy-Schwartz's inequality and Prop. \ref{propA2},
\begin{equation}
	\label{II5}
	\left\vert\varphi_q(D_x)p_k(x,D)(1-\psi(D))
	\varphi_p(D)u\right\vert_{L^1}\lesssim N^m_{\ind,0}(\sigma_{(q)})
	\left\vert
	\varphi_p(D)u\right\vert_{L^2}.
\end{equation}
From (\ref{II4}) and (\ref{II5}) one obtains
$$
	\left\vert R_p(p_k)u\right\vert_{L^1}\lesssim
	2^{-p(s+t)}\sum_{\vert p-q\vert\leq N}
	2^{qs}N^m_{\ind,0}(\sigma_{(q)})
	2^{pt}\left\vert 
	\varphi_p(D)u\right\vert_2,
$$
and the l.h.s. of (\ref{II3}) is therefore bounded from above by
\begin{eqnarray*}
	 \sum_{p\geq r-n_0} 2^{(r-p)(s+t)}
	\sum_{\vert p-q\vert\leq N}
	2^{qs}N^m_{\ind,0}(\sigma_{(q)})
	2^{pt}\left\vert  
	\varphi_p(D)u\right\vert_2.	
\end{eqnarray*}
Since $s+t>0$, H\"older's inequality yields
that the l.h.s. of (\ref{II2}) is bounded
from above by
$$
	\left\vert\left(
	2^{ps}N^m_{\ind,0}(\sigma_{(p)})
	2^{pt}\Big\vert  
	\varphi_p(D)u\right\vert_2\right)_{p\geq -1}\Big\vert_{l^1}.
$$
By Cauchy-Schwartz's inequality, Lemma \ref{lemmLP1} and an argument similar 
to the one used in (\ref{truc}), one obtains  (\ref{II2}), which concludes 
the proof.
\end{proof}
A first important consequence of Propositions \ref{propII3} and \ref{propII4}
is that one can control the action of the 
operator associated to the 'remainder'
symbol $\sigma-\slf-\sI$, which is  more regular than the full
operator if $\sigma(x,\xi)$ is smooth enough in the space variables.
\begin{prop}
	\label{propII5}
 	Let $m\in\R$, $s_0>d/2$ and $d/2<t_0\leq s_0$.
	If for some $r\geq 0$, one has $\sigma\in\Gamma^m_{s_0+r}$ 
	then,\\
	{\bf i.} For all $-t_0<s\leq s_0$, the following estimate
	holds:
	$$
	\forall u\in H^{m+t_0-r}(\R^d),\qquad
	\vert \Op(\sigma-\sI-\slf)u\vert_{H^s}\lesssim 
	N_{\ind,s+r}^m(\sigma)
	\vert u\vert_{H^{m+t_0-r}}.
	$$
	{\bf ii.} For all $r'\in\R$ (such that $t_0+r'\leq s_0+r$)
	and $-t_0<s\leq t_0+r'$, one has
	$$
	\forall u\in H^{s+m-r'}(\R^d),\qquad
	\vert \Op(\sigma-\sI-\slf)u\vert_{H^s}\lesssim 
	N_{\ind,t_0+r'}^m(\sigma)
	\vert u\vert_{H^{s+m-r'}}.
	$$
	{\bf iii.} For symbols of nonnegative order, i.e. when $m>0$, then
	for all $s>0$ such that $s+m\leq s_0$, one also has
	$$
	\forall u\in C_*^{-r}(\R^d),\qquad
	\vert \Op(\sigma-\sI-\slf)u\vert_{H^s}\lesssim
	N^m_{\ind,s+m+r}(\sigma)\vert u\vert_{C^{-r}_*};
	$$
	this estimate still holds for slightly negative values of $r$,
	namely, if $-m<r$.
\end{prop}
\begin{proof}
One has $\sigma-\sI-\slf=\sII+\sR$, and we 
are therefore led to control 
 $\vert \Op(\sII)u\vert_{H^s}$ and 
$\vert \Op(\sR)u\vert_{H^s}$. We first prove point {\bf i}.\\
The estimate on $\vert\Op(\sII)u\vert_{H^s}$ is given by the first
part of Prop. \ref{propII3} when $r=0$. When $r>0$, taking $s=s+r$ and $t=r$
in the second part of
this proposition gives the result.
The estimate on $\vert\Op(\sR)u\vert_{H^s}$ is given by taking $s=s+r$ and
$t=t_0-r$ in the first part of Prop. \ref{propII4}.\\
To establish {\bf ii.}, take $s=t_0+r'$ and $t=t_0-s+r'$ in the second 
estimate of Prop. \ref{propII3} to obtain that
$\vert \Op(\sII)u\vert_{H^s}\lesssim 
	N_{\ind,t_0+r'}^m(\sigma)
	\vert u\vert_{H^{s+m-r'}}$ for all $s<t_0+r'$. Taking $t=s-r'$ and
$s=t_0+r'$ in the first estimate of Prop. \ref{propII4} shows that
$\vert \Op(\sR)u\vert_{H^s}\lesssim 
	N_{\ind,t_0+r'}^m(\sigma)
	\vert u\vert_{H^{s+m-r'}}$ for all $s>-t_0$ and the proof of {\bf ii.}
is complete (the endpoint $s=t_0+r'$ being given by {\bf i}).\\
To prove {\bf iii.}, take $s=s+m+r$ and $t=m+r>0$
in the second part of Prop. \ref{propII3} and $s=s+m+r$ and $t=-m-r$ in
the second estimate of Prop. \ref{propII4}.
\end{proof}
The first two points of the following proposition are a close variant of 
Prop. \ref{propII5} which uses the decomposition (\ref{decompreste}) of
the component $\sR$, while the last point addresses the case when 
$\sigma(x,\xi)=\sigma(x)$ does not depend on $\xi$.
\begin{prop}
	\label{propII5bis}
 	Let $m\in\R$, $k\in\N$, $s_0>d/2$ and $d/2<t_0\leq s_0$.
	If for some $r\geq 0$, one has $\sigma\in\Gamma^m_{s_0+r}$ 
	then,\\
	{\bf i.} For all $-t_0<s\leq s_0$, the following estimate
	holds:
	$$
	\forall u\in H^{m+t_0-r},\quad 
	\vert \Op(\sigma-\sI-\slf-\sigma_{R,2})u\vert_{H^s}\lesssim 
	N_{\ind,s+r-k}^m(\nabla_x^k\sigma)
	\vert u\vert_{H^{m+t_0-r}}.
	$$
	{\bf ii.} For all $r'\in\R$ (such that $t_0+r'\leq s_0+r$)
	and $-t_0<s<t_0+r'$, one has
	$$
	\forall u\in H^{s+m-r'},\quad
	\vert \Op(\sigma-\sI-\slf-\sigma_{R,2})u\vert_{H^s}\lesssim 
	N_{\ind,t_0+r'-k}^m(\nabla_x^k\sigma)
	\vert u\vert_{H^{s+m-r'}}.
	$$
	{\bf iii.} For symbols of nonnegative order, i.e. when $m>0$, then
	for all $s>0$ such that $s+m\leq s_0$, one also has
	$$
	\forall u\in C_*^{-r}(\R^d),\qquad
	\vert \Op(\sigma-\sI-\slf-\sigma_{R,2})u\vert_{H^s}\lesssim
	N^m_{\ind,s+m+r-k}(\nabla^k_x\sigma)\vert u\vert_{C^{-r}_*}.
	$$
	{\bf iv.} When $\sigma$ is a function, $\sigma\in H^{s_0}(\R^d)$,
	one has
	$$
	\forall \quad 0<s\leq s_0, \qquad
	\vert \Op(\sigma-\sI-\slf-\sigma_{R,2})u\vert_{H^s}\lesssim
	\vert\nabla^k\sigma\vert_{H^{s-k}}\vert u\vert_{\infty}
	$$
	and, if $\sigma\in H^{s_0+r}(\R^d)$, with $r>0$,
	$$
	\forall \quad 0<s\leq s_0,\qquad
	\vert \Op(\sigma-\sI-\slf-\sigma_{R,2})u\vert_{H^s}\lesssim
	\vert\nabla^k\sigma\vert_{H^{s+r-k}}\vert u\vert_{C^{-r}_*};
	$$
	when $s=0$, the above two estimates still hold if one
	adds 
	$\vert\nabla^{n+1}\sigma\vert_{L^{\infty}}\vert u\vert_{H^{-n-1}}$
	to the right-hand-side, for any $n\in\N$.
\end{prop}
\begin{rem}
	\label{remII5bis}
	When $k=0$, the estimates of {\bf iv} still hold if one replaces
	$\Op(\sigma-\sI-\slf-\sigma_{R,2})$ by $\Op(\sigma-\sI-\slf)$
	(and $\vert \nabla^{n+1}\sigma\vert_\infty$ by
	$\vert \sigma\vert_{W^{n+1,\infty}}$ in the additional term
	when $s=0$).
	This is a consequence of the definition of $\sigma_{R,2}$ and
	of Lemma \ref{lemmII0}.
\end{rem}
\begin{proof}
One has $\sigma-\sI-\slf-\sigma_{R,2}=\sII+\sigma_{R,1}$, so that the
first three points of the proposition are proved as in Prop. \ref{propII5},
using Remarks \ref{remagrad} and \ref{remagrad1}.\\
We now prove the fourth point of the proposition. Since $\sigma$ is 
a function, we can write, as in Remark \ref{rempara}, 
$\Op(\sII+\sigma_{R,1})u=T_{\tu}\sigma+R(\widetilde{\sigma},\tu)$,
with $\tu:=(1-\psi(D))u$ and $\widetilde{\sigma}:=(1-\psi(D))\sigma$. 
The estimate for $s>0$ thus follows from the 
classical properties (e.g. Th. 2.4.1. of \cite{Chemin}, 
and Prop. 3.5.D of \cite{TaylorM0} for the last one):
\begin{itemize}
\item for all $s\in\R$,  
$\vert T_fg\vert_{H^s}\lesssim \vert f\vert_\infty \vert g\vert_{H^s}$;
\item for all $s\in\R$ and $r>0$,  
$\vert T_f g\vert_{H^s}\lesssim \vert f\vert_{C_*^{-r}}\vert g\vert_{H^{s+r}}$;
\item for all $s>0$, $r\in\R$,  
$\vert R(f,g)\vert_{H^s}\lesssim\vert f\vert_{H^{s+r}}\vert g\vert_{C_*^{-r}}$;
\item for all $n\in\N$,
$\vert R(f,g)\vert_{L^2}\lesssim
\vert  f\vert_{W^{n+1,\infty}}\vert g\vert_{H^{-n-1}}$;
\end{itemize}
(we also use the fact that $\vert f\vert_{H^s}\leq \cst \vert \nabla^k f\vert_{H^{s-k}}$ and $\vert f\vert_{W^{n,\infty}}\leq \cst \vert \nabla^n f\vert_{L^\infty}$ for all $f$ such that $\widehat{f}$ vanishes in a neighborhood of the origin).
\end{proof}
Gathering the results of the previous propositions, one obtains the following 
theorem, which describes the
action of the full operator $\Op(\sigma)$, which is of course of order $m$.
\begin{theo}
	\label{theoII1}
	Let $m\in\R$, $d/2<t_0\leq s_0$ and 
	$\sigma\in\Gamma_{s_0}^m$. Then for all $u\in\sch$,
	the following estimates hold:
	$$
	\forall -t_0<s<t_0,\qquad
	\vert \Op(\sigma)u\vert_{H^s}\lesssim \Big(
	n_{0,t_0}(\sigma)+N_{\ind,t_0}^{m}(\sigma)\Big)
	\vert u\vert_{H^{s+m}},
	$$
	and
	$$
	\forall t_0\leq s\leq s_0,\quad
	\vert \Op(\sigma)u\vert_{H^s}\lesssim \big(
	n_{0,s}(\sigma)+N_{\ind,s}^{m}(\sigma)\big)
	\vert u\vert_{H^{m+t_0}}+M_d^m(\sigma)\vert u\vert_{H^{s+m}}.
	$$
\end{theo}
\begin{proof}
Recall that $\sigma=\slf+\sI+(\sigma-\slf-\sI)$; we use the first two estimates
of Prop. \ref{propII5} (with $r=0$ and $r'=0$) to control
$\sigma-\slf-\sI$ while 
$\vert \Op(\slf)u\vert_{H^s}$ and 
$\vert \Op(\sI)u\vert_{H^s}$ are easily controlled using Propositions
\ref{propII1} and \ref{propII2} and the observation that by a classical
Sobolev embedding, $M_d^m(\sigma)\lesssim N_{\ind,t_0}^{m}(\sigma)$.
\end{proof}
\begin{rem}
	\label{remrem}
	{\bf i.} If $\sigma(x,\xi)=\sigma(x)\in H^{s_0}(\R^d)$, 
	then the results
	on the microlocal regularity of products (e.g. \cite{Hormander}, p.240)
	say that if $u\in H^s(\R^d)$, then $\sigma u\in H^s(\R^d)$ if
	$s+s_0>0$, $s\leq s_0$ and $s_0>d/2$. This result can be deduced from
	Th. \ref{theoII1} (note that the limiting case $s+s_0=0$ is also true,
	but the proof requires different tools \cite{Hormander}, th. 8.3.1).\\
	{\bf ii.} We refer to Prop. 8.1 of \cite{TaylorMpara} for another
	kind of estimate on the action of pseudo-differential operators;
	see also estimate (25) of \cite{Marschall}.
\end{rem}

One of the interests of Th. \ref{theoII1} is that it gives control of 
$\Op(\sigma)u$ in Sobolev spaces of negative order. The price to pay is
that for nonnegative values
of the Sobolev index $s$, and when 
$\sigma(\cdot,\xi)=\sigma(\cdot)\in H^{s_0}(\R^d)$ does not depend on $\xi$,
we do not recover the classical tame estimate 
$\vert \sigma u\vert_{H^s}\lesssim
\left(\vert u\vert_{H^s}\vert \sigma\vert_\infty
+\vert \sigma\vert_{H^s}\vert u\vert_\infty\right)$ but a weaker one, namely 
$\vert \sigma u\vert_{H^s}\lesssim \vert u\vert_{H^s}\vert \sigma\vert_{\infty}
+\vert \sigma\vert_{H^s}\vert u\vert_{H^{\frac{d}{2}+\varepsilon}}$, 
for all 
$\varepsilon>0$. The difference is slight because the embedding 
$H^{\frac{d}{2}+\varepsilon}(\R^d)\subset L^\infty(\R^d)$ is critical, but
can be cumbersome. The next theorem can therefore be a useful alternative to
Theorem \ref{theoII1}
\begin{theo}
	\label{theoII2}
	Let $m\in\R$, $d/2<t_0\leq s_0$ and 
	$\sigma\in\Gamma_{s_0+m}^m$. Assume that $m>0$ and $\sigma$ is 
	$\ind$-regular at the origin.
	Then for all $u\in\sch$ and $0<s\leq s_0$, one has
	$$
	\vert \Op(\sigma)u\vert_{H^s} \lesssim
	\Big(n_{\ind,s}(\sigma)+N^m_{\ind,s+m}(\sigma)\Big)
	\vert u\vert_{C_*^0}
	+M_d^m(\sigma)\vert u\vert_{H^{s+m}}.
	$$
\end{theo}
\begin{proof}
We just have to control
the $H^s$-norm of the four components of $\Op(\sigma)u$ by the r.h.s. of
the estimate given in the theorem.\\
For $\Op(\slf)u$ and $\Op(\sI)u$, this is a simple consequence of the
second part of Prop. \ref{propII1} and of Prop. \ref{propII2}
respectively. The other components are controlled with the help of 
Prop. \ref{propII5}.iii.
\end{proof}
\begin{rem}
	\label{exttheo}
	Using the fact that slightly negative values of $r$ are allowed
	in Prop. \ref{propII5}.iii, one can check 
	that the estimate given by Th. \ref{theoII2} 
	can be extended to $s=0$ provided that the quantity 
	$\vert u\vert_{C_*^0}$ which appears in the r.h.s. of the estimate 
	is replaced by $\vert u\vert_{C^\epsilon_*}$, for any $\epsilon>0$.
\end{rem}

The following corollary
deals with the case when the symbol $\sigma$ is of the form 
$\sigma(x,\xi)=\Sigma(v(x),\xi)$.
\begin{cor}
	\label{coroII1}
	Let $m\in\R$, $p\in\N$ and $s_0\geq t_0>d/2$. 
	Consider $v\in H^{s_0}(\R^d)^p$
	and assume that $\sigma(x,\xi)=\Sigma(v(x),\xi)$, with 
	$\Sigma\in C^\infty(\R^p,{\mathcal M}^m)$. Then:\\
	{\bf i.} 
	$$
	\forall -t_0< s <t_0,\qquad
	\vert \sigma(x,D)u\vert_{H^s}\lesssim
	C_\Sigma(\vert v\vert_\infty)\vert v\vert_{H^{t_0}}
	\vert u\vert_{H^{s+m}}
	$$
	and
	$$
	\forall t_0\leq s\leq s_0,\qquad
	\vert \sigma(x,D)u\vert_{H^s}\lesssim
	C_\Sigma(\vert v\vert_\infty)\left(\vert v\vert_{H^{s}}
	\vert u\vert_{H^{m+t_0}}
	+
	\vert u\vert_{H^{s+m}}\right).
	$$
	{\bf ii.} If moreover $m>0$, $v\in H^{s_0+m}(\R^d)$,
	 and $\Sigma$ is $\ind$-regular at the origin, then, 
	for $0<s\leq s_0$,
	$$
	\vert \sigma(x,D)u\vert_{H^s}\lesssim
	C_\Sigma(\vert v\vert_\infty)\big(\vert v\vert_{H^{s+m}}
	\vert u\vert_{C_*^0}
	+
	\vert u\vert_{H^{s+m}}\big).
	$$
	In the above, $C_\Sigma(\cdot)$ denotes a smooth nondecreasing function
	depending only on a finite number of derivatives of $\Sigma$.
\end{cor}
\begin{proof}
We write $\sigma(x,\xi)=[\sigma(x,\xi)-\Sigma(0,\xi)]+\Sigma(0,\xi)$.
Owing to Lemma \ref{lemmA0}, the first component of this decomposition 
is in $\Gamma^m_{s_0}$ and we can use Th. \ref{theoII1} to study the associated
pseudo-differential operator. The estimates of the theorem transform into
the estimates stated in the corollary thanks to Lemma \ref{lemmA0}.\\
Since the action of the Fourier multiplier $\Sigma(0,D)$ satisfies obviously
these estimates, the first point of the corollary is proved. Using Theorem
\ref{theoII2}, one proves the second estimate in the same way.
\end{proof}

%%%%%%%%%%%%%%%%%%%%%%%%%%%%%%%%%%%%%%%%%%%%%%%%%%%%%%%%%%%%%%%%%%%%%%%
%%%%%%%%%%%%%%%%%%%%%%%%%%%%%%%%%%%%%%%%%%%%%%%%%%%%%%%%%%%%%%%%%%%%%%%
%%%%%%%%%%%%%%%%%%%%%%%%%%%%%%%%%%%%%%%%%%%%%%%%%%%%%%%%%%%%%%%%%%%%%%%
\reseteq
\section{Composition and commutator estimates}\label{sectcomm}
%%%%%%%%%%%%%%%%%%%%%%%%%%%%%%%%%%%%%%%%%%%%%%%%%%%%%%%%%%%%%%%%%%%%%%%
%%%%%%%%%%%%%%%%%%%%%%%%%%%%%%%%%%%%%%%%%%%%%%%%%%%%%%%%%%%%%%%%%%%%%%%
%%%%%%%%%%%%%%%%%%%%%%%%%%%%%%%%%%%%%%%%%%%%%%%%%%%%%%%%%%%%%%%%%%%%%%%

The composition of two pseudo-differential operators is well-known for 
classical symbols, and one has 
$\Op(\sigma^1)\circ\Op(\sigma^2)\sim\Op(\sigma^1\sharp\sigma^2)$, where the 
symbol $\sigma^1\sharp\sigma^2$ is given by an infinite expansion of 
$\sigma^1$ and $\sigma^2$. When dealing
with symbols of limited regularity, one has to stop this expansion. 
Therefore, for all $n\in\N$, we define
$\sigma^1\sharp_n\sigma^2$ as
\begin{equation}
	\label{III1}
	\sigma^1\sharp_n\sigma^2(x,\xi):=\sum_{\vert\alpha\vert\leq n}
	\frac{(-i)^{\vert \alpha\vert}}{\alpha !}\partial_\xi^\alpha\sigma^1(x,\xi)
	\partial_x^\alpha\sigma^2(x,\xi).
\end{equation}
Similarly, we introduce the Poisson brackets:
\begin{equation}
	\label{IIIpoiss}
	\{\sigma^1,\sigma^2\}_n(x,\xi):=	
	\sigma^1\sharp_n\sigma^2(x,\xi)-
	\sigma^2\sharp_n\sigma^1(x,\xi).
\end{equation}

\bigbreak

In this section, we describe the composition or 
commutator of pseudodifferential operators
of limited regularity with Fourier multipliers, or with another
pseudo-differential operator.
A key point in this analysis is the following proposition; the first two
points are 
precise estimates for Meyer's classical result on the composition of
paradifferential operators (e.g. Th. XVI.4 of \cite{Meyer-Coifman3}).
\begin{prop}
	\label{lemmIII1}
	Let $m_1,m_2\in\R$, $s_0>d/2$, $n\in\N$ and 
	$\sigma^2\in \Gamma_{s_0+n+1}^{m_2}$. Then\\
	{\bf i.} If $\sigma^1(x,\xi)=\sigma^1(\xi)\in{\mathcal M}^{m_1}$,
	there exists a symbol $\rho_n(x,\xi)$ such that
	$$
	\sigma^1(D)\circ\Op(\sI^2)=
	\Op(\sigma^1\sharp_n\sI^2)+\Op(\rho_n);
	$$
	moreover $\rho_n(x,\xi)$ vanishes for $\vert\xi\vert\leq 1/2$ and
	satisfies the spectral condition (\ref{A3})
	and the estimate
	$$
	\sup_{\vert\xi\vert\geq 1/2}\sup_{\vert\beta\vert\leq d}\Big(
	\modj^{\vert\beta\vert+n+1-m_1-m_2}
	\vert\partial_\xi^\beta\rho_n(\cdot,\xi)\vert_\infty\Big)\lesssim 
	M^{m_1}_{n+2+[\frac{d}{2}]+d}(\sigma^1)
	M^{m_2}_{d}(\nabla_x^{n+1}\sigma^2).
	$$
	{\bf ii.} If $\sigma^1\in\Gamma_{s_0}^{m_1}$
	then there exists a symbol $\rho_n(x,\xi)$ such that
	$$
	\Op(\sI^1)\circ\Op(\sI^2)=\Op(\sI^1\sharp_n\sI^2)+\Op(\rho_n),
	$$
	and which satisfies the same properties as in 
	case {\bf i.}\\
	{\bf iii.} If $\sigma^1$ is a function, $\sigma^1\in C_*^r$ for
	some $r\geq 0$, then the symbol $\rho_n(x,\xi)$ defined in {\bf ii.}
	is of order
	$m_2-n-1-r$ and
	$$
	M_d^{m_2-n-1-r}(\rho_n)\lesssim 
	\vert \sigma^1\vert_{C_*^r}
	M^{m_2}_{d}(\nabla_x^{n+1}\sigma^2).
	$$
\end{prop}
\begin{rem}
	\label{remkey}
	For the sake of simplicity, we stated the above proposition for
	paradifferential symbols $\sI^1$ (and $\sI^2$ in ii. and iii.)
	 associated to 
	symbols $\sigma^1$ and $\sigma^2$; the proof below shows that
	the only specific properties of $\sI^1$ and $\sI^2$ actually used
	are the spectral property (\ref{A3}) and the cancellation for
	frequencies $\vert \xi\vert\leq 1/2$. Thus, one can extend the
	result to all symbols satisfying these conditions.
\end{rem}
\begin{proof}
We omit the proof of the first
point of the proposition, which can be be deduced from the proof below 
with only minor changes. The method we propose here is inspired by the proof
of Th. B.2.16 of \cite{Metivier-Zumbrun} rather than
Meyer's classical one (Th. XVI.4 of \cite{Meyer-Coifman3}) 
which would lead to less precise
estimates here.\\
First remark that since $\sI^1$ satisfies the spectral condition (\ref{A3})
and vanishes for frequencies $\vert\xi\vert\leq 1/2$,
there exists an admissible cut-off function $\chi$ (in the sense of
Def. \ref{defiLP0}) such that 
$\widehat{\sI^2}(\eta,\xi)\chi(\eta,\xi)=\widehat{\sI^2}(\eta,\xi)$; 
it is then both classical and easy to see  that 
$$
	\rho_n(x,\xi)=\sum_{\vert\gamma\vert=n+1}
	\int_{\R^d}G_\gamma(x,x-y,\xi)(\partial_x^\gamma\sI^2)(y,\xi)dy,
$$ 
with 
$$
	G_\gamma(x,y,\xi):=(-i)^{\vert\gamma\vert}(2\pi)^{-d}
	\int_{\R^d}e^{iy\cdot\eta}\sI^{1,\gamma}(x,\eta,\xi)
	\chi(\eta,\xi)d\eta
$$ 
and 
$\dsp \sI^{1,\gamma}(x,\eta,\xi):=
\int_0^1\frac{(1-t)^n}{n!}\partial_\xi^\gamma\sI^1(x,\xi+s\eta)ds$. 
Therefore, for all $0\leq\vert\beta\vert\leq d$,
\begin{eqnarray}
	\nonumber
	\vert \partial_\xi^\beta\rho_n(x,\xi)\vert&\leq&\cst
	\sum_{\beta'+\beta''=\beta}
	\left\vert \partial_\xi^{\beta'}G_\gamma(x,\cdot,\xi)\right\vert_{L^1}
	\left\vert \partial_x^\gamma\partial_\xi^{\beta''}\sI^2(\cdot,\xi)
	\right\vert_\infty\\
	\label{marre1}	
	&\leq&
	\cst 	\sum_{\beta'+\beta''=\beta}
	\left\vert\partial_\xi^{\beta'}G_\gamma(x,\cdot,\xi)\right\vert_{L^1}
	M_{d}^{m_2}(\nabla_x^{n+1}\sigma^2)\modj^{m_2-\vert\beta''\vert},
\end{eqnarray}
where we used Prop. \ref{propA1} to obtain the last equality.\\
The proposition follows therefore from (\ref{marre1}) and the estimate,
for all $\vert\beta'\vert\leq d$,
\begin{equation}
	\label{marre2}
	\left\vert\partial_\xi^{\beta'}G_\gamma(x,\cdot,\xi)\right\vert_{L^1}
	\leq  \cst M_{\Ind}^{m_1}(\sigma^1)
	\modj^{m_1-\vert\beta'\vert-n-1}.
\end{equation}
and, when $\sigma^1$ is a function (case {\bf iii.} of the lemma),
\begin{equation}
	\label{marre3}
	\left\vert\partial_\xi^{\beta'}G_\gamma(x,\cdot,\xi)\right\vert_{L^1}
	\leq\cst \vert \sigma^1\vert_{C_*^r}
	\modj^{-\vert\beta'\vert-n-1-r}.
\end{equation}
Both (\ref{marre2}) and (\ref{marre3}) follow from the next two lemmas.
\begin{lem}
	For all $\alpha,\beta\in\N^d$ such that
	$\vert\alpha\vert\leq [d/2]+1$ and $\vert\beta\vert\leq d$, one has
	$$
	\vert \partial_\eta^\alpha\partial_\xi^\beta
	(\sI^{1,\gamma}(x,\cdot,\cdot)\chi)(\eta,\xi)\vert
	\leq \cst M_{\Ind}^{m_1}(\sigma^1)
	\modj^{m_1-\vert\alpha\vert-\vert\beta\vert-n-1}.
	$$
	If moreover the symbol is a function, $\sigma^1\in C_*^r$ for
	some $r\in\R$, then
	$$
	\vert \partial_\eta^\alpha\partial_\xi^\beta
	(\sI^{1,\gamma}(x,\cdot,\cdot)\chi)(\eta,\xi)\vert
	\leq \cst \vert \sigma^1\vert_{C_*^r}
	\modj^{-\vert\alpha\vert-\vert\beta\vert-n-1-r}.
	$$
\end{lem}
\begin{proof}
It suffices to prove the estimate of the lemma for
$\partial_\eta^{\alpha'}\partial_\xi^{\beta'}\sI^{1,\gamma}
\partial_\eta^{\alpha''}\partial_\xi^{\beta''}\chi$ for all 
$\alpha'+\alpha''=\alpha$ and $\beta'+\beta''=\beta$. By definition of
$\sI^{1,\gamma}$, one has
$$
	\partial_\eta^{\alpha'}\partial_\xi^{\beta'}\sI^{1,\gamma}
	\partial_\eta^{\alpha''}\partial_\xi^{\beta''}\chi(\eta,\xi)
	=
	\int_0^1\frac{(1-t)^n}{n!}\partial_\xi^{\alpha'+\beta'+\gamma}
	\sI^1(x,\xi+s\eta)s^{\vert\alpha'\vert}ds \,	
	\partial_\eta^{\alpha''}\partial_\xi^{\beta''}\chi(\eta,\xi).
$$
Since on the support of $\partial_\eta^{\alpha''}\partial_\xi^{\beta''}\chi$
one has $\langle \xi+s\eta\rangle\sim \modj$,
the first estimate of the lemma follows from the definition of the 
seminorms $M_k(\cdot)$ and (\ref{estadm}).\\
When $\sigma^1$ is a function, and since 
$\vert\alpha'\vert+\vert\beta'\vert+\vert\gamma\vert\geq 1$, we can use 
Prop. A.5 of \cite{Auscher-Taylor} which asserts that
$\vert\partial_\xi^{\alpha'+\beta'+\gamma}\sI^1(x,\xi)\vert
\leq \cst \vert \sigma^1\vert_{C_*^r}\modj^{-r-\vert\alpha'\vert
-\vert\beta'\vert-\vert\gamma\vert}$, from which one easily obtains the 
second estimate of the lemma.
\end{proof}
\begin{lem}
	Let $F(\cdot,\cdot)$ be a function defined on $\R^d_\eta\times\R^d_\xi$
	and such that
	\begin{itemize}
	\item There exists $0<\delta<1$ such that $F(\eta,\xi)=0$
	for all $\vert\eta\vert\geq \delta\vert\xi\vert$;
	\item 	For all $\alpha\in\N^d$, $\vert\alpha\vert\leq[d/2]+1$,
	there exists a constant $C_\alpha$
	such that
	$$
	\forall \eta,\xi\in\R^d,\qquad
	\left\vert \partial_\eta^\alpha F(\eta,\xi)\right\vert
	\leq C_\alpha \modj^{\mu-\vert\alpha\vert}.
	$$
	\end{itemize}
	Then, one has
	$$
	\forall\xi\in\R^d,\qquad
	\left\vert \fourinv{F}(\cdot,\xi)\right\vert_{L^1}
	\leq \cst \Big(\sup_{\vert\alpha\vert\leq [d/2]+1}C_\alpha\Big)
	\modj^\mu.
	$$
\end{lem}
\begin{proof}
This result can be proved with the techniques used to prove the estimate
(2.21) of Appendix B in \cite{Metivier-Zumbrun}. Briefly, and for the sake
of completeness, we sketch the proof. Define 
$F^\flat(\eta,\xi):=F(\modj \eta,\xi)$ and remark that 
$\vert \fourinv{F}(\cdot,\xi)\vert_{L^1}=
\vert \fourinv{F^\flat}(\cdot,\xi)\vert_{L^1}\leq \cst 
\vert F^\flat(\cdot,\xi)\vert_{H^{[d/2]+1}}$. The first assumption made in the
statement of the lemma shows that $F^\flat(\cdot,\xi)$ is supported in the
ball $\{\vert\eta\vert\leq 1\}$ so that it is easy to conclude using the
second assumption.
\end{proof}
\end{proof}

%%%%%%%%%%%%%%%%%%%%%%%%%%%%%%%%%%%%%%%%%%%%%%%%%%%%%%%%%%%%%%%%%%%%%%%
%%%%%%%%%%%%%%%%%%%%%%%%%%%%%%%%%%%%%%%%%%%%%%%%%%%%%%%%%%%%%%%%%%%%%%%
\subsection{Commutators with Fourier multipliers}\label{sectfour}
%%%%%%%%%%%%%%%%%%%%%%%%%%%%%%%%%%%%%%%%%%%%%%%%%%%%%%%%%%%%%%%%%%%%%%%
%%%%%%%%%%%%%%%%%%%%%%%%%%%%%%%%%%%%%%%%%%%%%%%%%%%%%%%%%%%%%%%%%%%%%%%

We give in this section some commutator estimates between a 
Fourier multiplier and a pseudo-differential 
operator of limited regularity. We first set some notations:

\bigbreak

\noindent
{\bf Notations.} For all $m\in\R$, $s_0>d/2$, and all symbols 
$\sigma\in\Gamma^m_{s_0}$,
we define
\begin{equation}
	\label{defnorm1}
	\forall s\leq s_0,\qquad
	\Vert \sigma \Vert_{H^s_{(m)}}:=
	\frac{1}{2}\left( n_{0,s}(\sigma)+N_{\ind,s}^m(\sigma)\right)
\end{equation}
and, when $\sigma$ is also $\ind$-regular at the 
origin,
\begin{equation}
	\label{defnorm2}
	\forall s\leq s_0,\qquad
	\Vert \sigma \Vert_{H^s_{reg,(m)}}:=
	\frac{1}{2}\left( n_{\ind,s}(\sigma)+N_{\ind,s}^m(\sigma)\right).
\end{equation}
Finally, if $\sigma$ is $d$-regular at the origin, we set
\begin{equation}
	\label{defnormnoel}
	\Vert \sigma\Vert_{\infty,(m)}:=\frac{1}{2}
	\left(
	m_d(\sigma)+M^m_d(\sigma)
	\right).
\end{equation}
When no confusion is possible, we omit the subscript $m$ in 
the above definitions. 
Remark that
when $\sigma$ does not depend on $\xi$ (i.e. when it is a function), then
one has $\Vert\sigma\Vert_{H^s}=\Vert\sigma\Vert_{H^s_{reg}}
=\vert \sigma\vert_{H^s}$, and 
$\Vert\sigma\Vert_\infty=\vert \sigma\vert_\infty$.

\bigbreak

The first commutator estimates we state are of Kato-Ponce type:
\begin{theo}
	\label{theoIII1}
	Let $m_1,m_2\in\R$, $n\in\N$ and $d/2<t_0\leq s_0$.
	Let $\sigma^1(\xi)\in{\mathcal M}^{m_1}$ be $n$-regular
	at the origin and 
	$\sigma^2(x,\xi)\in\Gamma^{m_2}_{s_0+m_1\wedge n+1}$. Then:\\
	{\bf i.} For all $s\in\R$ such that
	$\max\{-t_0,-t_0-m_1\}<s\leq s_0+1$, one has
	\begin{eqnarray*}
	\lefteqn{\big\vert \big[\sigma^1(D),\Op(\sigma^2)\big]u
	-\Op(\{\sigma^1,\sigma^2\}_n)u
	\big\vert_{H^s}}\\
	&\lesssim & C(\sigma^1)\big(
	M_{d}^{m_2}(\nabla_x^{n+1}\sigma^2)
	\vert u\vert_{H^{s+m_1+m_2-n-1}}
	+\Vert \sigma^2\Vert_{H^{s+m_1\wedge n}}
	\vert u\vert_{H^{m_1+m_2+t_0-m_1\wedge n}}\big),
	\end{eqnarray*}
	where $C(\sigma^1):=
	M^{m_1}_{n+2+[\frac{d}{2}]+d}(\sigma^1)+m_n(\sigma^1)$.\\
	{\bf ii.} If moreover $\sigma^1$ is $(\Ind)$-regular
	 and $\sigma^2$ is $d$-regular at the origin,
	then the above estimate can be replaced by
	\begin{eqnarray*}
	\lefteqn{\big\vert \big[\sigma^1(D),\Op(\sigma^2)\big]u
	-\Op(\{\sigma^1,\sigma^2\}_n)u
	\big\vert_{H^s}}\\
	&\lesssim & C'(\sigma^1)\big(
	\Vert \nabla_x^{n+1}\sigma^2\Vert_\infty
	\vert u\vert_{H^{s+m_1+m_2-n-1}}
	+\Vert \nabla_x^{n+1}\sigma^2\Vert_{H^{s+m_1\wedge n-n-1}}
	\vert u\vert_{H^{m_1+m_2+t_0-m_1\wedge n}}\big),
	\end{eqnarray*}
	where $C'(\sigma^1):=
	M^{m_1}_{n+2+[\frac{d}{2}]+d}(\sigma^1)+m_{\Ind}(\sigma^1)$.
\end{theo}
\begin{proof}
{\bf i.} Remark that ${\{\sigma^1,\sigma^2\}_n}
=\sigma^1{\sharp_n}\sigma^2-\sigma^1\sigma^2$ 
and that $\Op(\sigma^2)\circ\Op(\sigma^1)=\Op(\sigma^1\sigma^2)$ since 
$\Op(\sigma^1)$ is a Fourier multiplier. Therefore, one has 
$$
	\left[\Op(\sigma^1),\Op(\sigma^2)\right]-\Op(\{\sigma^1,\sigma^2\}_n)
	=\Op(\sigma^1)\circ\Op(\sigma^2)-\Op(\sigma^1\sharp_n\sigma^2).
$$
Write now $\sigma^1(D)\circ\Op(\sigma^2)-\Op(\sigma^1{\sharp_n}\sigma^2)=
\sum_{j=1}^5\tau^j(x,D)$, with
\begin{eqnarray*}
	\tau^1(x,D)&=&\sigma^1(D)\circ\Op(\sigma^2-\sI^2-\slf^2),\\
	\tau^2(x,D)&= &\sigma^1(D)\circ \Op(\slf^2),\\
	\tau^3(x,D)&=&\sigma^1(D)\circ \Op(\sI^2)
	-\Op(\sigma^1\sharp_n\sigma_I^2),\\
	\tau^4(x,D)&=&\Op(\sigma^1\sharp_n\sI^2
	-(1-\psi(\xi))\sigma^1{\sharp_n}\sigma^2),\\
	\tau^5(x,D)&=&-\Op(\psi(\xi)\sigma^1\sharp_n\sigma^2).
\end{eqnarray*}
We now turn to control the operator norms of $\tau^j(x,D)$, $j=1,\dots,5$.\\
$\bullet$ {\bf Control of $\mathbf{\tau^1(x,D)}$.} Since $\sigma^1(D)$ is a
Fourier multiplier, one obtains easily that for all $s\in\R$ and 
$u\in\sch$, 
$$
	\vert \tau^1(x,D)u\vert_{H^s}\leq (m_0(\sigma^1)+M_0^{m_1}(\sigma^1))
	\left\vert \Op(\sigma^2-\sI^2-\slf^2)u\right\vert_{H^{s+m_1}}.
$$
Using Prop. \ref{propII5}.i (with $r=m_1\wedge n-m_1$) gives therefore,
for all $-t_0-m_1<s\leq s_0+1$,
\begin{equation}
	\label{bor1}
	 \vert \tau^1(x,D)u\vert_{H^s}\lesssim 
	(m_0(\sigma^1)+M_0^{m_1}(\sigma^1))
	N^{m_2}_{\ind,s+m_1\wedge n}(\sigma^2)
	\vert u\vert_{H^{m_1+m_2+t_0-m_1\wedge n}}.
\end{equation}
$\bullet$ {\bf Control of $\mathbf{\tau^2(x,D)}$.} One has
$$
	\vert \tau^2(x,D)u\vert_{H^s}\leq 
	(m_0(\sigma^1)+M^{m_1}_0(\sigma^1))
	\vert \Op(\slf^2)u\vert_{H^{s+m_1}},
$$ 
so that it
is a direct consequence of Prop. \ref{propII1} that one has, for all 
$s\leq s_0+1$,
\begin{equation}
	\label{bor2}
	\vert \Op(\tau^2)u\vert_{H^s}\lesssim 
	(m_0(\sigma^1)+M_0^{m_1}(\sigma^1))
	n_{0,s+m_1}(\sigma^2) \vert u\vert_{H^{m_1+m_2+t_0-m_1\wedge n}}.
\end{equation}
$\bullet$ {\bf Control of $\mathbf{\tau^3(x,D)}$.} We have 
$\Op(\tau^3)=\Op(\rho_n)$ with $\rho_n$ as given in the first part of 
Prop. \ref{lemmIII1}. This lemma asserts that the symbol $\rho_n(x,\xi)$ 
satisfies the conditions of application of Lemma \ref{lemmII1}, which states 
that $\vert \Op(\rho_n)u\vert_{H^s}\lesssim
M^{m_1+m_2-n-1}_d(\rho_n)\vert u\vert_{H^{s+m_1+m_2-n-1}}$. 
Using the estimate of $M^{m_1+m_2-n-1}_d(\rho_n)$ given in 
Prop. \ref{lemmIII1} shows therefore that for all $s\in\R$,
\begin{equation}
	\label{bor3}
	\vert \Op(\tau^3)u\vert_{H^s}\lesssim
	M^{m_1}_{\Ind}(\sigma^1)M^{m_2}_{d}(\nabla_x^{n+1}\sigma^2)
	\vert u\vert_{H^{s+m_1+m_2-n-1}}.
\end{equation}
$\bullet$ {\bf Control of $\mathbf{\tau^4(x,D)}$.} By definition of the product 
law $\sharp_n$, one has
\begin{eqnarray*}
	\tau^4(x,\xi)&=&
	-\sum_{\vert\alpha\vert\leq n}\frac{(-i)^\alpha}{\alpha!}
	\partial_\xi^\alpha \sigma^1(\xi)
	\partial_x^\alpha (\sigma^2-\sI^2-\slf^2)(x,\xi)\\
	&=&-\sum_{\vert\alpha\vert\leq n}\frac{(-i)^\alpha}{\alpha!}
	{\mathbf 1}_{[1/2,\infty)}(\xi)
	\partial_\xi^\alpha \sigma^1(\xi)
	\partial_x^\alpha (\sigma^2-\sI^2-\slf^2)(x,\xi),
\end{eqnarray*}
where ${\mathbf 1}_{[1/2,\infty)}(\cdot)$ denotes the characteristic function
of the interval $[1/2,\infty)$. It follows that
$$
	\big\vert \tau^4(x,D)u\big\vert_{H^s}\lesssim
	\sum_{\vert\alpha\vert\leq n}
	\big\vert\Op(\partial_x^\alpha\sigma^2-(\partial_x^\alpha\sigma^2)_I
	-(\partial_x^\alpha\sigma^2)_{lf})v
	\big\vert_{H^s}
$$
with $v=\Op({\mathbf 1}_{[1/2,\infty)}\partial_\xi^\alpha\sigma^1)u$; 
we now use
the first estimate of Prop. \ref{propII5} 
(with $r=m_1\wedge n-\vert\alpha\vert$) to obtain that the terms of the 
above sum are bounded from
above by $N^{m_2}_{\ind,s+m_1\wedge n-\vert\alpha\vert}
(\partial_x^\alpha\sigma^2)\vert v
\vert_{H^{m_2+t_0-m_1\wedge n+\vert\alpha\vert}}$, for all
$-t_0<s\leq s_0+1$. It is now straightforward to conclude that
for all $-t_0<s\leq s_0+1$,
\begin{equation}
	\label{bor5}
	\big\vert \tau^4(x,D)u\big\vert_{H^s}\lesssim
	M^{m_1}_{n}(\sigma^1) N^{m_2}_{\ind,s+m_1\wedge n}(\sigma^2)
	\vert u\vert_{H^{m_1+m_2+t_0-m_1\wedge n}}.
\end{equation}
$\bullet$ {\bf Control of $\mathbf{\tau^5(x,D)}$.} One has 
$\Op(\psi(\xi)\sigma^1\sharp_n\sigma^2)u
=\Op((\sigma^1\sharp_n\sigma^2)_{lf})u$, so that
Prop. \ref{propII1} can be used to obtain for all $s\leq s_0+1$,
\begin{equation}
	\label{bor6}
	\left\vert\Op(\psi(\xi)\sigma^1\sharp_n\sigma^2)u\right\vert_{H^s}
	\lesssim m_n(\sigma^1)n_{0,s+n}(\sigma^2)
	\vert u\vert_{H^{m_1+m_2+t_0-m_1\wedge n}}.
\end{equation}

Recalling that $\Vert \sigma^2\Vert_{H^s}$ is defined in (\ref{defnorm1}), 
the estimate given in {\bf i.} of the theorem now follows directly from 
(\ref{bor1}), (\ref{bor2}), (\ref{bor3}),
(\ref{bor5}) and (\ref{bor6}).\\
{\bf ii.} We use here another decomposition, namely 
$\sigma^1(D)\circ\Op(\sigma^2)-\Op(\sigma^1{\sharp_n}\sigma^2)=
\sum_{j=1}^6\ut^j(x,D)$, with
\begin{eqnarray*}
	\ut^1(x,D)&=&\sigma^1(D)\circ
	\Op(\sigma^2-\sI^2-\slf^2-\sigma^2_{R,2}),\\
	\ut^2(x,D)&= &\sigma^1(D)\circ \Op(\Psi(D_x)\slf^2+\sigma^2_{R,2})
	-\Op(\sigma^1\sharp_n(\psi(D_x)\sigma^2_{lf}+\sigma^2_{R,2})),\\
	\ut^3(x,D)&=&\sigma^1(D)\circ \Op(\sI^2)
	-\Op(\sigma^1\sharp_n\sigma_I^2),\\
	\ut^4(x,D)&=&-\Op(\sigma^1\sharp_n(\sigma^2-
	\sI^2-\slf^2
	-\sigma^2_{R,2})),\\
	\ut^5(x,D)&=&\sigma^1(D)\circ\Op((1-\psi(D_x)\sigma^2_{lf})),\\
	\ut^6(x,D)&=&-\Op(\sigma^1\sharp_n(1-\psi(D_x))\slf^2).
\end{eqnarray*}
We now turn to control the operator norms of $\ut^j(x,D)$, $j=1,\dots,6$.\\
{\bf Control of $\mathbf{\ut^1(x,D)}$.} Proceeding as for the control of 
$\tau^1(x,D)$ in {\bf i.} above, but using Prop. \ref{propII5bis} 
instead of Prop. \ref{propII5}, one can replace 
$N^{m_2}_{\ind,s+m_1\wedge n}(\sigma^2)$ by 
$N^{m_2}_{\ind,s+m_1\wedge n-n-1}(\nabla_x^{n+1}\sigma^2)$ in (\ref{bor1}).\\
{\bf Control of $\mathbf{\ut^2(x,D)}$.} We need here two lemmas:
\begin{lem}
	Let $m_1\in\R$, $n\in\N$
	and $\sigma^1(\xi)\in{\mathcal M}^{m_1}$ be $(\Ind)$-regular
	at the origin. Let 
	$\sigma^2(x,\xi)$ be a symbol $d$-regular at the origin and
	such that $\widehat{\sigma^2}(\eta,\xi)$ is supported in
	the ball $\vert \eta\vert+\vert\xi\vert\leq A$, for some $A>0$.\\
	Then, $\sigma^1(D)\circ \Op(\sigma^2)=\Op(\sigma^1\sharp_n\sigma^2)
	+\Op(\rho_n)$, where the symbol $\rho_n(x,\xi)$ is such that
	$\widehat{\rho_n}(\eta,\xi)$ vanishes outside the ball 
	$\vert\eta\vert+\vert\xi\vert\leq A$ and satisfies the estimate
	\begin{eqnarray*}
	\sup_{\vert\xi\vert\leq A}\sup_{\vert\beta\vert\leq d} 
	\vert \partial_\xi^\beta\rho_n(\cdot,\xi)\vert_\infty
	&\leq& \cst
		\sup_{\vert\xi\vert\leq 2A}
	\sup_{\vert\alpha\vert\leq n+2+[\frac{d}{2}]+d}
	\vert\partial_\xi^\alpha\sigma^1(\xi)\vert\\
	& &\times\sup_{\vert\xi\vert\leq A}\sup_{\vert\alpha\vert\leq d}
	\vert\nabla^{n+1}_x\partial_\xi^\alpha
	\sigma^2(\cdot,\xi)\vert_\infty.
	\end{eqnarray*}
\end{lem}
\begin{proof}
The proof is a close adaptation of the proof of Prop. \ref{lemmIII1}. First
replace the admissible cut-off function $\chi(\eta,\xi)$ used there by
a smooth function $\widetilde{\chi}(\eta,\xi)$ supported in the ball
$\vert\eta\vert+\vert\xi\vert\leq A$. Inequality (\ref{marre1}) must then
be replaced by
$$
	\vert\partial_\xi^\beta\rho_n(x,\xi)\vert\leq\cst
	\Big(\sum_{\vert\beta'\vert\leq d}\big\vert\partial_\xi^{\beta'}G_\gamma(\cdot,\xi)\big\vert_{L^1}\Big)
	\sup_{\vert\xi\vert\leq A}\sup_{\vert \alpha\vert\leq d}\vert\nabla^{n+1}_x\partial_\xi^\alpha\sigma^2(\cdot,\xi)\vert_\infty,
$$
for all $\vert\beta\vert\leq d$ and $\vert\xi\vert\leq A$.\\
Finally, one concludes the proof as in Prop. \ref{lemmIII1} after remarking that (\ref{marre2}) can be replaced here by
$$
	\big\vert\partial_\xi^{\beta'}G_\gamma(\cdot,\xi)\big\vert_{L^1}
	\leq\cst
	\sup_{\vert\xi\vert\leq 2A}
	\sup_{\vert\alpha\vert\leq n+2+[\frac{d}{2}]+d}
	\vert\partial_\xi^\alpha\sigma^1(\xi)\vert.
$$
\end{proof}
The proof of the following lemma is a straightforward adaptation of the proof
of Lemma \ref{lemmII1}.
\begin{lem}
	Let $\rho(x,\xi)$ be a symbol such that $\widehat{\rho}(\eta,\xi)$ is
	supported in the ball $\vert\eta\vert+\vert\xi\vert\leq A$, for some
	$A>0$. Then $\Op(\rho)$ extends as a continuous operator on every
	Sobolev space, with values in $H^\infty(\R^d)$. Moreover,
	$$
	\forall s,t\in\R,\forall u\in H^t(\R^d),\qquad
	\vert \Op(\rho)u\vert_{H^s}\leq \cst
	\sup_{\vert\xi\vert\leq A}\sup_{\vert\alpha\vert\leq d} 
	\vert \partial_\xi^\alpha\rho(\cdot,\xi)\vert_\infty
	\vert u\vert_{H^t},
	$$
	where the constant depends only on $A$, $s$ and $t$.
\end{lem}
From these two lemmas, one easily gets
$$
	\vert \ut^2(x,D)u\vert_{H^s}\lesssim C'(\sigma^1)
	\Vert \nabla_x^{n+1}\sigma^2\Vert_\infty
	\vert u\vert_{H^{s+m_1+m_2-n-1}}.
$$
where $C'(\sigma^1)$ is as given in the statement of the theorem.\\
{\bf Control of $\mathbf{\ut^3(x,D)}$.} One has $\ut^3=\tau^3$, so that
$\ut^3(x,D)$ is controlled via (\ref{bor3}).\\
{\bf Control of $\mathbf{\ut^4(x,D)}$.} To control this term, proceed exactly as
for the control of $\tau^4(x,D)$ above, but use Prop. \ref{propII5bis}
instead of Prop. \ref{propII5}. This yields
$$
	\big\vert \ut^4(x,D)u\big\vert_{H^s}\lesssim
	M^{m_1}_{n}(\sigma^1) 
	N^{m_2}_{\ind,s+m_1\wedge n-n-1}(\nabla_x^{n+1}\sigma^2)
	\vert u\vert_{H^{m_1+m_2+t_0-m_1\wedge n}}.
$$
{\bf Control of $\mathbf{\ut^5(x,D)}$ and  $\mathbf{\ut^6(x,D)}$.} 
The difference between $\ut^5(x,D)$
and $\tau^2(x,D)$ is that the operator $(1-\psi(D_x))$ is applied to 
$\slf^2(\cdot,\xi)$ in the former. This allows us to replace $n_{0,s+m_1}(\sigma^2)$   by	$n_{0,s+m_1-n-1}(\nabla_x^{n+1}\sigma^2)$ 
in (\ref{bor2}).\\
A similar adaptation of (\ref{bor6}) gives 
$$
	\big\vert\Op(\ut^6)u\big\vert_{H^s}
	\lesssim m_n(\sigma^1)n_{0,s-1}(\nabla_x^{n+1}\sigma^2)
	\vert u\vert_{H^{m_1+m_2+t_0-m_1\wedge n}}.
$$

Point {\bf ii} of the theorem thus follows from the estimates
on $\ut^j(x,D)$, $j=1,\dots,6$, proved above.
\end{proof}
\begin{rem}
	When $m_2\geq 0$, an easy adaptation of the above proof shows that
	the quantity 
	$\Vert \sigma^2\Vert_{H^{s+m_1\wedge n}}
	\vert u\vert_{H^{m_1+m_2+t_0-m_1\wedge n}}$ 
	which appears in the r.h.s. of the first estimate 
	of the theorem can 
	be replaced by 
	$\Vert \sigma^2\Vert_{H^{s+m_1\wedge n+m_2}}
	\vert u\vert_{H^{t_0+m_1-m_1\wedge n}}$.
\end{rem}

In the spirit of Theorem \ref{theoII2}, the following two theorems can be a
useful alternative to Theorem \ref{theoIII1}. Theorem \ref{theoIII2} deals
with the case of pseudo-differential operators $\sigma^2(x,D)$ of nonnegative
order, while Th. \ref{theoIII3} addresses the case 
$\sigma^2(x,\xi)=\sigma^2(x)$.
\begin{theo}
	\label{theoIII2}
	Let $m_1\in\R$, $m_2>0$, $n\in\N$ and $s_0>d/2$. Let
	$\sigma^1(\xi)\in{\mathcal M}^{m_1}$ be $n$-regular
	at the origin and let 
	$\sigma^2(x,\xi)\in \Gamma^{m_2}_{s_0+m_1\wedge n+1}$
	be $\ind$-regular at the origin. Then:\\
	{\bf i.} For all $s$ such that
	$0<s+m_1$, $0<s$ and $s+m_2\leq s_0+1$, one has
	\begin{eqnarray*}
	\lefteqn{\big\vert \big[\sigma^1(D),\Op(\sigma^2)\big]u
	-\Op(\{\sigma^1,\sigma^2\}_n)u
	\big\vert_{H^s}}\\
	&\lesssim& C(\sigma^1)\big(M_{d}^{m_2}(\nabla_x^{n+1}\sigma^2)
	\vert u\vert_{H^{s+m_1+m_2-n-1}}
	+
	\Vert \sigma^2\Vert_{H^{s+m_1\wedge n+m_2}_{reg}}
	\vert u\vert_{\infty}\big),
	\end{eqnarray*}
	where $C(\sigma^1)=M^{m_1}_{n+2+[\frac{d}{2}]+d}(\sigma^1)
	+m_n(\sigma^1)$.\\
	{\bf ii.} If moreover $\sigma^1$ is $(\Ind)$-regular
	at the origin,
	then the above estimate can be replaced by
	\begin{eqnarray*}
	\lefteqn{\big\vert \big[\sigma^1(D),\Op(\sigma^2)\big]u
	-\Op(\{\sigma^1,\sigma^2\}_n)u
	\big\vert_{H^s}}\\
	&\lesssim& C'(\sigma^1)\big(\Vert \nabla_x^{n+1}\sigma^2\Vert_\infty
	\vert u\vert_{H^{s+m_1+m_2-n-1}}
	+
	\Vert \nabla_x^{n+1}\sigma^2\Vert_{H^{s+m_1\wedge n+m_2-n-1}_{reg}}
	\vert u\vert_{\infty}\big),
	\end{eqnarray*}
	where $C'(\sigma^1):=
	M^{m_1}_{n+2+[\frac{d}{2}]+d}(\sigma^1)+m_{\Ind}(\sigma^1)$.
\end{theo}
\begin{proof}
We only give a sketch of the proof of {\bf i.} which follows the same lines 
as the 
proof of Theorem \ref{theoIII1}.i; {\bf ii.} is obtained similarly. 
The modifications to be made are:\\
- Inequality (\ref{bor1}) must be replaced by
\begin{equation}
	\label{modif1}
	 \vert \tau^1(x,D)u\vert_{H^s}\lesssim 
	(m_0(\sigma^1)+M_0^{m_1}(\sigma^1))
	N^{m_2}_{\ind,s+m_1+m_2}(\sigma^2)\vert u\vert_{\infty},
\end{equation}
which holds for all $s\in\R$ such that $s+m_1>0$ and $s+m_2\leq s_0+1$.
This is a consequence of Prop. \ref{propII5}.iii,
which can be used since we assumed $m_2>0$.\\
- Similarly, the second estimate of Prop. \ref{propII1} allows us to
replace (\ref{bor2}) by
$$
	\vert \Op(\tau^2)u\vert_{H^s}\lesssim 
	(m_0(\sigma^1)+M_0^{m_1}(\sigma^1))
	n_{\ind,s+m_1}(\sigma^2) \vert u\vert_{\infty}.
$$
- Inequality (\ref{bor3}) is left unchanged.\\
- Estimate (\ref{bor5}) must be replaced
by
\begin{equation}
	\label{modif2}
	\big\vert \tau^4(x,D)u\big\vert_{H^s}\lesssim
	M^{m_1}_n(\sigma^1) N^{m_2}_{\ind,s+m_2+m_1\wedge n}(\sigma^2)
	\vert u\vert_{\infty},
\end{equation}
which holds for all $0<s$ and $s+m_2\leq s_0+1$. One proves (\ref{modif2}) in the same
way as (\ref{bor5}), using the third point of Prop. \ref{propII5} rather
than the first one (this is possible because $m_2>0$ here).\\
- Finally, inequality (\ref{bor6}) is replaced, using the second part of 
Prop. \ref{propII1}, by
$$
	\big\vert\Op(\psi(\xi)\sigma^1\sharp_n\sigma^2)u\big\vert_{H^s}
	\lesssim m_n(\sigma^1)n_{\ind,s+n}(\sigma^2)
	\vert u\vert_{\infty}.
$$
\end{proof}
An interesting particular case is obtained when the symbol $\sigma^2(x,\xi)$
does not depend on $\xi$ (i.e., it is a function).
\begin{theo}
	\label{theoIII3}
	Let $m_1\in\R$, $n\in\N$ and $s_0>d/2$. Let
	$\sigma^1(\xi)\in{\mathcal M}^{m_1}$ be $n$-regular
	at the origin and let 
	$\sigma^2\in H^{s_0+m_1\wedge n+1}(\R^d)$.\\
	{\bf i.} If $m_1>n$ then for all $s$ such that
	$0\leq s\leq s_0+1$, one has
	\begin{eqnarray*}
	\lefteqn{\big\vert \big[\sigma^1(D),\sigma^2\big]u
	-\Op(\{\sigma^1,\sigma^2\}_n)u
	\big\vert_{H^s}}\\
	&\leq &C(\sigma^1)\big(
	\vert\sigma^2\vert_{W^{n+1,\infty}}
	\vert u\vert_{H^{s+m_1-n-1}}
	+\vert \sigma^2\vert_{H^{s+m_1}}
	\vert u\vert_{\infty}\big),
	\end{eqnarray*}
	where $C(\sigma^1)=
	M^{m_1}_{n+2+[\frac{d}{2}]+d}(\sigma^1)+m_n(\sigma^1)$.\\
	{\bf ii.} If moreover $\sigma^1$ is $(\Ind)$-regular at the origin,
	then the above estimate can be replaced by
	\begin{eqnarray*}
	\lefteqn{\big\vert \big[\sigma^1(D),\sigma^2\big]u
	-\Op(\{\sigma^1,\sigma^2\}_n)u
	\big\vert_{H^s}}\\
	&\leq &C'(\sigma^1)\big(
	\vert\nabla^{n+1}\sigma^2\vert_{\infty}
	\vert u\vert_{H^{s+m_1-n-1}}
	+\vert \nabla^{n+1}\sigma^2\vert_{H^{s+m_1-n-1}}
	\vert u\vert_{\infty}\big),
	\end{eqnarray*}
	where $C'(\sigma^1):=
	M^{m_1}_{n+2+[\frac{d}{2}]+d}(\sigma^1)+m_{\Ind}(\sigma^1)$.
\end{theo}
\begin{proof}
Here again, we only prove the first point of the theorem since the proof 
of the second one can be deduced similarly from the proof 
of Th. \ref{theoIII1}.ii.\\
Remark that for all
$k\in\N$, $s\in\R$ and $v\in\sch$, one has 
$N^0_{k,s}(v)=\vert v\vert_{H^s}$. Therefore, we just have to adapt the 
points of the proof of Th. \ref{theoIII2} which
use the assumption $m_2>0$, namely the obtention of (\ref{modif1}) and 
(\ref{modif2}).\\
One can check that (\ref{modif1}) remains true here. This is a consequence
of Prop. \ref{propII5bis}.iv and Remark \ref{remII5bis}.\\
We now prove that (\ref{modif2}) can be replaced by 
$$
 	\vert \tau^4(x,D)u\vert_{H^s}\lesssim M_n^{m_1}(\sigma^1)
	\vert \sigma^2\vert_{H^{s+m_1}}\vert u\vert_{\infty},
$$
which holds for all $0<s\leq s_0+1$ and provided that $m_1>n$, and which
remains true when $s=0$ provided one adds 
$\vert \sigma^2\vert_{W^{n+1,\infty}}\vert u\vert_{H^{m_1-n-1}}$ 
to the right-hand-side.
To obtain this, we need to control in $H^s$-norm, and for all 
$0\leq \vert\alpha\vert\leq n$, the terms
$\Op(\partial_x^\alpha\sigma^2-(\partial_x^\alpha \sigma^2)_{I}
-(\partial_x^\alpha\sigma^2)_{lf})(\partial_\xi^\alpha \sigma^1(D))u$, 
which is done using Prop. \ref{propII5bis}.iv (with $r=m_1-\vert\alpha\vert>0$)
and Remark \ref{remII5bis}.
\end{proof}

We finally give commutator estimates of Calderon-Coifman-Meyer type:
\begin{theo}
	\label{theoIII1bis}
	Let $m_1,m_2\in\R$, $n\in\N$ and $d/2<t_0\leq s_0$.
	Let $\sigma^1(\xi)\in{\mathcal M}^{m_1}$ be $n$-regular
	at the origin and 
	$\sigma^2(x,\xi)\in\Gamma^{m_2}_{s_0+m_1\wedge n+1}$. Then:\\
	{\bf i.} For all $s\in\R$ such that	
	$-t_0<s\leq t_0+1$ and $-t_0<s+m_1\leq t_0+n+1$,
	$$
	\big\vert \big[\sigma^1(D),\Op(\sigma^2)\big]u
	-\Op(\{\sigma^1,\sigma^2\}_n)u
	\big\vert_{H^s}
	\lesssim  C(\sigma^1)
	\Vert \sigma^2\Vert_{H^{t_0+n+1}}
	\vert u\vert_{H^{s+m_1+m_2-n-1}},
	$$
	where $C(\sigma^1):=
	M^{m_1}_{n+2+[\frac{d}{2}]+d}(\sigma^1)+m_n(\sigma^1)$.\\
	{\bf ii.} If moreover $\sigma^1$ is $(\Ind)$-regular
	 and $\sigma^2$ is $d$-regular at the origin,
	then the above estimate can be replaced by
	$$
	\big\vert \big[\sigma^1(D),\Op(\sigma^2)\big]u
	-\Op(\{\sigma^1,\sigma^2\}_n)u
	\big\vert_{H^s}
	\lesssim  C'(\sigma^1)
	\Vert \nabla_x^{n+1}\sigma^2\Vert_{H^{t_0}}
	\vert u\vert_{H^{s+m_1+m_2-n-1}},
	$$
	where $C'(\sigma^1):=
	M^{m_1}_{n+2+[\frac{d}{2}]+d}(\sigma^1)+m_{\Ind}(\sigma^1)$.
\end{theo}
\begin{proof}
The proof follows closely the proof of Th. \ref{theoIII1}, so that we
just mention the adaptations that have to be made.\\
{\bf i.} Below is the list of changes one must perform in the 
control of the operators  $\tau^j(x,D)$.\\
$\bullet$ {\bf Control of $\mathbf{\tau^1(x,D)}$.}  For all 
$-t_0<s+m_1\leq t_0+n+1$ and using 
Prop. \ref{propII5}.ii with $r'=n+1$ instead of 
Prop. \ref{propII5}.i, one obtains instead of (\ref{bor1}) a control
in terms of $N^{m_2}_{\ind,t_0+n+1}(\sigma^2)
	\vert u\vert_{H^{s+m_1+m_2-n-1}}$.\\
$\bullet$ {\bf Control of $\mathbf{\tau^2(x,D)}$.} For $s+m_1\leq
 t_0+n+1$, one just has to remark that Prop. \ref{propII1} gives a 
control in terms of 
$n_{0,t_0+n+1}(\sigma^2) \vert u\vert_{H^{s+m_1+m_2-n-1}}$.\\
$\bullet$ {\bf Control of $\mathbf{\tau^4(x,D)}$.} 
If $-t_0<s\leq t_0+1$, one can replace (\ref{bor5}) by a control in terms of
$N^{m_2}_{\ind,t_0+n+1}(\sigma^2)
	\vert u\vert_{H^{s+m_1+m_2-n-1}}$ provided that one
uses Prop. \ref{propII5}.ii with $r'=n+1-\vert\alpha\vert$ instead
of Prop. \ref{propII5}.i.\\
$\bullet$ {\bf Control of $\mathbf{\tau^5(x,D)}$.} When $s\leq t_0+1$ one gets
easily a control in terms of $\Vert \sigma^2\Vert_{H^{t_0+1}}\vert u\vert_{H^{s+m_1+m_2-n-1}}$.\\
{\bf ii.} One deduces the second point of the theorem from Th. \ref{theoIII1}.ii
exactly as we adapted the proof of the first point 
from the proof of Th. \ref{theoIII1}.i.
\end{proof}

\begin{rem}
	{\bf i.} Extending results of Moser \cite{Moser} and 
	Kato-Ponce \cite{Kato-Ponce}, Taylor proved in \cite{TaylorM0} the 
	following generalized Kato-Ponce 
	estimates (which also hold in $L^p$-based Sobolev spaces): 
	for all Fourier multiplier $\sigma^1(D)$ of order
	$m_1>0$, and all $\sigma^2\in H^\infty(\R^d)$, one has
	for all $s\geq0$,
	\begin{equation}
	\label{estTay}
	\big\vert [\Op(\sigma^1),\sigma^2]u\big\vert_{H^s}
	\leq \underline{C}(\sigma^1)\big(\vert\nabla\sigma^2\vert_{\infty} 
	\vert u\vert_{H^{s+m_1-1}}
	+\vert \sigma^2\vert_{H^{s+m_1}}\vert u\vert_\infty\big),
	\end{equation}
	where $\underline{C}(\sigma^1)$ is some 
	constant depending on $\sigma^1$
	(in \cite{TaylorM0}, Taylor also deals with classical 
	pseudo-differential operators $\sigma^1(x,D)$ --and not 
	only Fourier multipliers--; we address
	this problem in Section \ref{sectfin}).\\
	The estimate of Th. \ref{theoIII3} coincides with (\ref{estTay})
	when $n=0$ (it is in fact more precise since it allows one to
	replace the term $\vert \sigma^2\vert_{H^{s+m_1}}$ by 
	$\vert \nabla\sigma^2\vert_{H^{s+m_1-1}}$ \footnote{Such an improvement
	has been proved recently in \cite{Benzoni-Danchin-Descombes} 
	for $n=0$ or $n=1$, 
	in the case where $\sigma^1(\xi)=\modj^{m_1}$ and 
	$\sigma^2$ does not depend on $\xi$}); 
	the general case $n\in \N$ gives an extension of this
	result involving the Poisson bracket of $\sigma^1$ and $\sigma^2$. 
	Th. \ref{theoIII2}
	extends (\ref{estTay}) in another
	way, allowing $\sigma^2(x,D)$ to be a pseudo-differential operator
	of nonnegative order instead of a function. Finally, the most
	general extension of (\ref{estTay}) is given by Th. \ref{theoIII1},
	since it contains the two improvements just mentioned, allows
	estimates in 
	Sobolev spaces of negative order, and does not assume 
	cumbersome restrictions
	on the order of $\sigma^1(\xi)$ and $\sigma^2(x,\xi)$. For instance,
	(\ref{estTay}) does not hold when $s<0$ or $m_1\leq 0$ but 
	can be replaced by:
	for all Fourier multiplier  $\sigma^1(D)$ of order
	$m_1\in\R$ regular at the origin, 
	and all $\sigma^2\in H^\infty(\R^d)$, one has for 
	all $t_0>d/2$ and $s>max\{-t_0,-t_0-m_1\}$,
	$$
	\big\vert [\Op(\sigma^1),\sigma^2]u\big\vert_{H^s}
	\leq C(\sigma^1)\big(\vert \nabla\sigma^2\vert_\infty 
	\vert u\vert_{H^{s+m_1-1}}
	+\vert \nabla\sigma^2\vert_{H^{s+(m_1)_+-1}}
	\vert u\vert_{H^{t_0+(m_1)_-}}\big),
	$$
	with $(m_1)_+=\max\{m_1,0\}$ and  $(m_1)_-=\min\{m_1,0\}$.\\
	{\bf ii.} The Calderon-Coifman-Meyer commutator estimate
	of Th. \ref{theoIII1bis} coincides with (\ref{intro3}) when
	$n=0$ and $\sigma^1$ is a Fourier multiplier (the general case
	is addressed in Section \ref{sectfin}) but its range of validity is
	wider since it allows negative values of $s$ and $m_1$. We also
	have the same kind of generalization as for the Kato-Ponce estimates.
\end{rem}

In the particular case when the symbol $\sigma^2$ is of the form 
$\sigma^2(x,\xi)=\Sigma(v(x),\xi)$, one can now obtain easily, proceeding
as in the proof of Corollary \ref{coroII1}:
\begin{cor}
	\label{coroIII3}
	Let $m_1,m_2\in\R$, $p\in\N$ and $d/2<t_0\leq s_0$.
	Let $\sigma^1(\xi)\in{\mathcal M}^{m_1}$ be $n$-regular
	at the origin
	and 
	$\sigma^2(x,\xi)=\Sigma^2(v(x),\xi)$ with
	$\Sigma^2\in C^\infty(\R^p,{\mathcal M}^{m_2})$ and 
	$v\in H^{s_0+m_1\wedge n+1}(\R^d)^p$. Then,\\
	{\bf i.} For all $s\in\R$ such that 
	$\max\{-t_0,-t_0-m_1\}<s\leq s_0+1$ 
	\begin{eqnarray*}
	\lefteqn{\big\vert \big[\sigma^1(D),\sigma^2(x,D)\big]u
	-\Op({\{\sigma^1,\sigma^2\}_n})
	u\big\vert_{H^s}}\\
	&\leq& C(\sigma^1)C_{\Sigma^2}(\vert v\vert_{W^{n+1,\infty}})
	\big(\vert u\vert_{H^{s+m_1+m_2-n-1}}
	+
	\vert v\vert_{H^{(s+m_1\wedge n)_+}}
	\vert u\vert_{H^{m_1+m_2+t_0-m_1\wedge n}}\big).
	\end{eqnarray*}
	{\bf ii.} For all $s\in\R$ such that 
	$-t_0<s<t_0+1$ and $-t_0<s+m_1\leq t_0+n+1$, one also has
	$$
	\big\vert \big[\sigma^1(D),\sigma^2(x,D)\big]u
	-\Op({\{\sigma^1,\sigma^2\}_n})
	u\big\vert_{H^s}
	 \leq  C(\sigma^1)C_{\Sigma^2}(\vert v\vert_{\infty})
	\vert v\vert_{H^{t_0+n+1}}\vert u\vert_{H^{s+m_1+m_2-n-1}}.
	$$
	{\bf iii.} If moreover $m_2>0$ and $\Sigma^2$ is $\ind$-regular
	at the origin, 
	one has, for all $s\in\R$ such that $0<s+m_1$, $0<s$, and 
	$s+m_2\leq s_0+1$,
	\begin{eqnarray*}
	\lefteqn{\big\vert \big[\sigma^1(D),\sigma^2(x,D)\big]u
	-\Op({\{\sigma^1,\sigma^2\}_n})
	u\big\vert_{H^s}}\\
	& \leq& C(\sigma^1)C_{\Sigma^2}(\vert v\vert_{W^{n+1,\infty}})
	\big(\vert u\vert_{H^{s+m_1+m_2-n-1}}
	+
	\vert v\vert_{H^{s+m_2+m_1\wedge n}}\vert u\vert_{\infty}\big).
	\end{eqnarray*}
	In the above, $C_{\Sigma^2}(\cdot)$ denotes a 
	smooth nondecreasing function
	depending only on a finite  number 
	of derivatives of $\Sigma^2$ and 
	$C(\sigma^1)=M^{m_1}_{n+2+[\frac{d}{2}]+d}(\sigma^1)+m_n(\sigma^1)$.
\end{cor}

\begin{exam}
	Let $\sigma(x,\xi)$ be the symbol given by (\ref{symbDN}),
	with $a\in H^\infty(\R^d)$, 
	and let $m\geq 0$.
	Then for all $t_0>d/2$ and $s>-t_0$, one has
	$$
	\left\vert \left[\langle D\rangle^m,\sigma(x,D)\right]u
	\right\vert_{H^s}
	\leq C(\vert \nabla a\vert_{W^{1,\infty}})
	\left(	
	\vert u\vert_{H^{s+m}}+\vert \nabla a\vert_{H^{(s+m)_+}}
	\vert u\vert_{H^{1+t_0}}\right),
	$$
	and, writing $\widetilde{m}:=\max\{m,1\}$
	$$
	\left\vert \left[\langle D\rangle^m,\Op (\sigma)\right]u
	-\Op(\tau)u\right\vert_{H^s}
	\leq C(\vert \nabla a\vert_{W^{2,\infty}})
	\left(	
	\vert u\vert_{H^{s+m-1}}+\vert \nabla a\vert_{H^{(s+\widetilde{m})_+}}
	\vert u\vert_{H^{1+t_0}}\right),
	$$
	where the symbol $\tau(x,\xi)$ is given by
	$$
	\tau(x,\xi):=m\langle \xi\rangle^{m-2}
	\frac{\vert\xi\vert^2 d^2a(\nabla a,\xi)-
	(\nabla a\cdot \xi)d^2a(\xi,\xi)}{\sigma(x,\xi)}.
	$$
	For $s>0$, we also have
	$$
	\left\vert \left[\langle D\rangle^m,\sigma(x,D)\right]u
	\right\vert_{H^s}
	\leq C(\vert \nabla a\vert_{W^{1,\infty}})
	\left(	
	\vert u\vert_{H^{s+m}}+\vert \nabla a\vert_{H^{1+s+m}}
	\vert u\vert_{{\infty}}\right).
	$$
	Finally, if $s$ and $t_0$ are such that $-t_0< s$ and 
	$s+m\leq t_0+1$, then
	$$
	\left\vert \left[\langle D\rangle^m,\sigma(x,D)\right]u
	\right\vert_{H^s}
	\leq C(\vert \nabla a\vert_{W^{1,\infty}})
	\vert \nabla a\vert_{H^{t_0+1}}\vert u\vert_{H^{s+m}}.
	$$
\end{exam}

%%%%%%%%%%%%%%%%%%%%%%%%%%%%%%%%%%%%%%%%%%%%%%%%%%%%%%%%%%%%%%%%%%%%%%%
%%%%%%%%%%%%%%%%%%%%%%%%%%%%%%%%%%%%%%%%%%%%%%%%%%%%%%%%%%%%%%%%%%%%%%%
\subsection{Composition and commutators between two 
pseudo-differential operators of limited regularity}\label{sectfin}
%%%%%%%%%%%%%%%%%%%%%%%%%%%%%%%%%%%%%%%%%%%%%%%%%%%%%%%%%%%%%%%%%%%%%%%
%%%%%%%%%%%%%%%%%%%%%%%%%%%%%%%%%%%%%%%%%%%%%%%%%%%%%%%%%%%%%%%%%%%%%%%
This section is devoted to the proof of composition and 
commutator estimates involving
two pseudo-differential operators of limited regularity. We first
introduce the following notations:

\bigbreak

\noindent
{\bf Notations.} For all $s_0>d/2$, $n\in\N$, and all symbols 
$\sigma\in\Gamma^m_{s_0}$ $n$-regular at the origin,
we define
\begin{equation}
	\label{defnorm3}
	\forall s\leq s_0,\qquad
	\Vert \sigma \Vert_{H^s_{n,(m)}}:=
	\frac{1}{2}\left( n_{n,s}(\sigma)+N_{\Ind,s}^m(\sigma)\right)
\end{equation}
and
\begin{equation}
	\label{defnorm4}
	\forall k\in\N, \quad 0\leq k\leq n,\qquad
	\Vert \sigma \Vert_{W^{k,\infty}_{n,(m)}}:=
	\frac{1}{2}\left( m_{n}(\sigma)+M_{\Ind,k}^m(\sigma)\right).
\end{equation}
When no confusion is possible, we omit the subscript $m$. Remark that
when $\sigma$ does not depend on $\xi$ (i.e. when it is a function), then
one has $\Vert\sigma\Vert_{H^s_n}=\vert \sigma\vert_{H^s}$ and 
$\Vert \sigma\Vert_{W^{k,\infty}_n}=\vert \sigma\vert_{W^{k,\infty}}$. 

\bigbreak

We first give commutator estimates of Kato-Ponce type:
\begin{theo}
	\label{theoIV1}
	Let $m_1,m_2\in\R$, $n\in\N$ and $d/2<t_0\leq s_0$.
	Define $m:=m_1\wedge m_2$ and
	let $\sigma^j\in \Gamma^{m_j}_{s_0+m\wedge n+1}$ ($j=1,2$) be two 
	symbols $n$-regular
	at the origin.\\
	{\bf i.} For all  $s$ such that
	$\max\{-t_0,-t_0-m_1\}<s\leq t_0+1$, one has
	\begin{eqnarray*}
	\lefteqn{\big\vert \Op(\sigma^1)\circ \Op(\sigma^2)u
	-\Op(\sigma^1\sharp_n\sigma^2)u
	\big\vert_{H^s}}\\
	&\lesssim&\Vert\sigma^1\Vert_{W^{n+1,\infty}_n}
	\Vert\sigma^2\Vert_{W^{n+1,\infty}}
	\vert u\vert_{H^{s+m_1+m_2-n-1}}\\
	&+&\big(\Vert \sigma^1\Vert_{H^{t_0+1}_n}
	\Vert\sigma^2\Vert_{H^{s_++m_1\wedge n}}
	+\Vert \sigma^1\Vert_{H^{s_++m_1\wedge n}_n}	
	\Vert \sigma^2\Vert_{H^{t_0}}\big)
	\vert u\vert_{H^{m_1+m_2+t_0-m_1\wedge n}},
	\end{eqnarray*}
	where $s_+=\max\{s,0\}$ and with the notations 
	(\ref{defnorm1})-(\ref{defnorm2}) and 
	(\ref{defnorm3})-(\ref{defnorm4}).\\
	{\bf i.bis.} Under the same assumptions, one also has
	\begin{eqnarray*}
	\lefteqn{\big\vert \Op(\sigma^1)\circ \Op(\sigma^2)u
	-\Op(\sigma^1\sharp_n\sigma^2)u
	\big\vert_{H^s}}\\
	&\lesssim&\Vert \sigma^1\Vert_{H^{t_0+m_1\wedge n+1}_n}\big(
	\Vert\sigma^2\Vert_{H^{s+m_1\wedge n}}
	\vert u\vert_{H^{m_1+m_2+t_0-m_1\wedge n}}
	+\Vert\sigma^2\Vert_{W^{n+1,\infty}}
	\vert u\vert_{H^{s+m_1+m_2-n-1}}\big).
	\end{eqnarray*}
	{\bf i.ter.} For $s\in\R$ such that 
	$max\{t_0+1,-t_0-m_1\}\leq s \leq s_0+1$,
	the estimate of {\bf i.} still holds if one adds
	$\Vert \sigma^1\Vert_{H^s}\Vert \sigma^2\Vert_{H^{m_1\wedge n+t_0+1}}		\vert u\vert_{H^{m_1+m_2+t_0-m_1\wedge n}}$
	to the right-hand--side.\\
	{\bf ii.} For all $s\in\R$ such that
	$\max\{-t_0,-t_0-m_1,-t_0-m_2\}<s\leq t_0+1$, one has
	$$
	\big\vert \big[\Op(\sigma^1),\Op(\sigma^2)\big]u
	-\Op(\{\sigma^1,\sigma^2\}_n)u
	\big\vert_{H^s}\lesssim I(\sigma^1,\sigma^2,m_1,m_2)+
	I(\sigma^2,\sigma^1,m_2,m_1),
	$$
	where $I(\sigma^1,\sigma^2,m_1,m_2)$ denotes the r.h.s. or
	the estimate of {\bf i}, {\bf i.bis} or {\bf i.ter.} In this
	latter case, the range of application of the estimate is bounded
	from above by $s_0+1$ instead of $t_0+1$.\\
	{\bf iii.} If $m_1>0$, $\sigma^1$ is $\ind$ regular at the origin,		and  $\sigma^2(x,\xi)=\sigma^2(x)$ does not 
	depend on $\xi$ then, for all $0\leq s\leq t_0+1$,
	$$
	\big\vert \big[\Op(\sigma^1),\Op(\sigma^2)\big]u
	\big\vert_{H^s}
	\lesssim\Vert \sigma^1\Vert_{H^{t_0+1}_{reg}}\big(
	\vert\sigma^2\vert_{H^{s+m_1}}
	\vert u\vert_{\infty}
	+\vert\sigma^2\vert_{W^{1,\infty}}
	\vert u\vert_{H^{s+m_1-1}}\big).
	$$
\end{theo}
\begin{proof}
One has 
\begin{eqnarray*}
	[\Op(\sigma^1),\Op(\sigma^2)]-\Op(\{\sigma^1,\sigma^2\})&=&
	\Op(\sigma^1)\circ\Op(\sigma^2)-\Op(\sigma^1\sharp_n\sigma^2)\\
	& &-\left(\Op(\sigma^2)\circ\Op(\sigma^1)
	-\Op(\sigma^2\sharp_n\sigma^1)\right)\\
	&:=&\tau_1(x,D)-\tau_2(x,D).
\end{eqnarray*}
Controlling the operator norm of $\tau_1(x,D)$ gives the composition 
estimates of
the theorem. To obtain the commutators estimates, we need also a 
control on $\tau_2(x,D)$; since this latter is obtained by a simple 
permutation of 
$\sigma^1$ and $\sigma^2$, we just have to treat $\tau_1(x,D)$. We
 decompose this operator into 
$ \tau_1(x,D)=\sum_{j=1}^7\tau_1^j(x,D)$
with
\begin{eqnarray*}
	\tau_1^1(x,D)&:=&\Op(\sigma^1)\circ\Op(\sigma^2-\slf^2-\sI^2)\\
	\tau_1^2(x,D)&:=&\Op(\sigma^1)\circ\Op(\slf^2)\\
	\tau_1^3(x,D)&:= &
	\left[\Op(\sI^1)\circ\Op(\sI^2)-\Op(\sI^1\sharp_n\sI^2)\right]\\
	\tau_1^4(x,D)&:= &\left[\Op(\sI^1\sharp_n\sI^2)
	-\Op((1-\psi(\xi))\sigma^1\sharp_n\sigma^2)\right]\\
	\tau_1^5(x,D)&:= &\Op((\sigma^1\sharp_n\sigma^2)_{lf})\\
	\tau_1^6(x,D)&:=&\Op(\slf^1)\circ\Op(\sI^2)\\
	\tau_1^7(x,D)&:=&\Op(\sigma^1-\slf^1-\sI^1)\circ \Op(\sI^2).
\end{eqnarray*}
The proof reduces therefore to
the control of $\vert \tau_1^j(x,D)u\vert_{H^s}$ for all $j=1,\dots,7$.\\
$\bullet$ {\bf Control of $\mathbf{\tau_1^1(x,D)}$ and $\mathbf{\tau_1^2(x,D)}$.} Using the first estimate of 
Th. \ref{theoII1}, one gets that for all $-t_0<s\leq t_0+1$,
\begin{equation}
	\label{blob}
	\vert\tau_1^j(x,D)u\vert_{H^s}
	\lesssim \Vert \sigma^1\Vert_{H^{t_0+1}}
	\vert v_j\vert_{H^{s+m_1}},	
	\qquad (j=1,2)
\end{equation}
with $v_1:=\Op(\sigma^2-\slf^2-\sI^2)u$ and $v_2:=\Op(\slf^2)u$.\\
Proceeding as for the estimates of $\tau^1$ and $\tau^2$ in the proof
of Th. \ref{theoIII1}, one obtains, for all 
$\max\{-t_0,-t_0-m_1\}<s\leq t_0+1$,
\begin{equation}
	\label{IV1}
	\vert \tau_1^j(x,D)u\vert_{H^s}\lesssim
	\Vert \sigma^1\Vert_{H^{t_0+1}}\Vert \sigma^2\Vert_{H^{s+m_1\wedge n}}
	\vert u\vert_{H^{m_1+m_2+t_0-m_1\wedge n}}.
\end{equation}
$\bullet$ {\bf Control of $\mathbf{\tau_1^3(x,D)}$.} A direct use of 
Prop. \ref{lemmIII1} and Lemma \ref{lemmII1} yields, for all $s\in\R$,
\begin{equation}
	\label{IV3}
	\vert\tau_1^3(x,D)u\vert_{H^s}\lesssim
	M^{m_1}_{\Ind}(\sigma^1)M_{d}^{m_2}(\nabla_x^{n+1}\sigma^2)
	\vert u\vert_{H^{s+m_1+m_2-n-1}}.
\end{equation}
$\bullet$ {\bf Control of $\mathbf{\tau_1^4(x,D)}$.} Obviously, it suffices
to control the operator norm of 
$A_\alpha(x,D):=\Op(\partial_\xi^\alpha\sigma^1_I\partial_x^\alpha\sigma_I^2)
-\Op((1-\psi(\xi))\partial_\xi^\alpha\sigma^1\partial_x^\alpha\sigma^2)$ 
for all $0\leq\vert\alpha\vert\leq n$. Let us introduce here $\tN:=N+3$; we
can assume that the paradifferential decomposition (\ref{A6}) used in this
proof is done using the integer $\tN$ instead of $N$. To enhance this
fact, we write momentaneously $\stI$ the paradifferential symbol associated
to any symbol $\sigma$ using (\ref{A7}) with $N$ replaced by $\tN$. We
denote $\sI$ when $N$ is used. We can write 
$A_\alpha(x,D)=A_\alpha^1(x,D)+A_\alpha^2(x,D)$ with
\begin{eqnarray}
	\label{compo1}
	A_\alpha^1(x,D)&:=&
	\Op\left(\partial_\xi^\alpha\stI^1\partial_x^\alpha\stI^2
	-(\partial_\xi^\alpha\sigma^1\partial_x^\alpha\sigma^2)_I\right),\\
	\label{compo2}
	A_\alpha^2(x,D)&:=&-
	\Op\left(\partial_\xi^\alpha\sigma^1\partial_x^\alpha\sigma^2
	-(\partial_\xi^\alpha\sigma^1\partial_x^\alpha\sigma^2)_{lf}
	-(\partial_\xi^\alpha\sigma^1\partial_x^\alpha\sigma^2)_I\right).
\end{eqnarray}
The operator norm of $A_\alpha^1(x,D)$ is controlled using the next lemma,
whose proof is postponed to Appendix \ref{applemm} to ease the readability
of the present proof. Note that this result is in the spirit of the main
result of \cite{Yamazaki2}, but that the estimate given in this latter 
reference is not useful here.
\begin{lem}
	\label{lemmfond}
	Let $\sigma^1(x,\xi)$ and $\sigma^2(x,\xi)$ be as in the 
	statement of the theorem, and let $\alpha\in\N^d$ be
	such that $0\leq\vert\alpha\vert\leq n$. Then,
	$$
	\forall u\in \sch,\qquad \left\vert A^1_\alpha(x,D)u
	\right\vert_{H^s}
	\lesssim
	M^{m_1}_{n+d,n+1}(\sigma^1)M^{m_2}_{d,n+1}(\sigma^2)
	\vert u\vert_{H^{s+m_1+m_2-n-1}},
	$$
	where $A_\alpha^1(x,D)$ is defined in (\ref{compo1}).
\end{lem}
To control $A_\alpha^2(x,D)$, we use Prop. \ref{propII5}.i. with 
$r=n\wedge m_1-\vert \alpha\vert$ to obtain, for all $-t_0\leq s\leq s_0+1$,
$$
	\left\vert A^2_\alpha(x,D)u
	\right\vert_{H^s}
	\lesssim
	N^{m_1+m_2-\vert\alpha\vert}_{\ind,s+m_1\wedge n-\vert\alpha\vert}
	(\partial_\xi^\alpha\sigma^1\partial_x^\alpha\sigma^2)
	\vert u\vert_{H^{m_1+m_2+t_0-m_1\wedge n}}.
$$
Using classical tame product estimates, one gets easily
\begin{equation}
	\label{adapt}
	N^{m_1+m_2-\vert\alpha\vert}_{\ind,s+m_1\wedge n-\vert\alpha\vert}
	(\partial_\xi^\alpha\sigma^1\partial_x^\alpha\sigma^2)
	\lesssim
	\Vert \sigma^1\Vert_{W^{0,\infty}_n}
	\Vert \sigma^2\Vert_{H^{s_+ +m_1\wedge n}}
	+
	\Vert \sigma^1\Vert_{H^{s_+ +m_1\wedge  n}_n}
	\Vert \sigma^2\Vert_{W^{0,\infty}_0}.
\end{equation}
We have thus proved the following estimate on $\Op(\tau^4_1)$, for
all $-t_0<s\leq s_0+1$:
\begin{eqnarray}
	\label{IV4}
	\vert \Op(\tau^4_1)u\vert_{H^s}\lesssim
	M^{m_1}_{n+d,n+1}(\sigma^1)M^{m_2}_{d,n+1}(\sigma^2)
	\vert u\vert_{H^{s+m_1+m_2-n-1}}\\
	\nonumber
	+\big(\Vert \sigma^1\Vert_{W^{0,\infty}_n}
	\Vert \sigma^2\Vert_{H^{s_+ +m_1\wedge n}}
	+
	\Vert \sigma^1\Vert_{H^{s_+ +m_1\wedge  n}_n}
	\Vert \sigma^2\Vert_{W^{0,\infty}_0}\big)\vert u\vert_{H^{m_1+m_2+t_0-m_1\wedge n}} 
\end{eqnarray}
$\bullet$ {\bf Control of $\mathbf{\tau_1^5(x,D)}$.} Using again classical
tame product estimates, one gets 
$$
	n_{0,s}((\sigma^1\sharp_n\sigma^2)_{lf})
	\lesssim
	m_n(\sigma^1)n_{0,s_++n}(\sigma^2)
	+n_{n,s_++n}(\sigma^1)m_{0,0}(\sigma^2)
$$
so that by Prop. \ref{propII1}, one obtains, for all $s\leq s_0+1$,
\begin{equation}
	\label{IV5}
	\vert \Op(\tau^5_1)u\vert_{H^s}\lesssim
	\big(m_n(\sigma^1)n_{0,s_++n}(\sigma^2)
	+n_{n,s_++n}(\sigma^1)m_{0,0}(\sigma^2)\big)
	\vert u\vert_{H^{m_1+m_2+t_0-m_1\wedge n}}.
\end{equation}
$\bullet$ {\bf Control of $\mathbf{\tau_1^6(x,D)}$.} Using successively 
Props. \ref{propII1} and \ref{propII2}
 one obtains, for all $s\leq s_0+n+1$,
\begin{equation}
	\label{IV6}
	\vert \Op(\tau^6_1)u\vert_{H^s}\lesssim
	n_{0,s}(\sigma^1)M_d^{m_2}(\sigma^2)
	\vert u\vert_{H^{m_1+m_2+t_0-m_1\wedge n}}.
\end{equation}
$\bullet$ {\bf Control of $\mathbf{\tau_1^7(x,D)}$.} Using successively 
Props. \ref{propII5}.i (with $r=m_1\wedge n$) and \ref{propII2}
 one obtains, for all $-t_0<s\leq s_0+n+1$,
\begin{equation}
	\label{IV7}
	\vert \Op(\tau^7_1)u\vert_{H^s}\lesssim
	N_{0,s+m_1\wedge n}^{m_1}(\sigma^1)M_d^{m_2}(\sigma^2)
	\vert u\vert_{H^{m_1+m_2+t_0-m_1\wedge n}}.
\end{equation}
$\bullet$ {\bf Proof of i.} Gathering the estimates
(\ref{IV1})-(\ref{IV3}) and (\ref{IV4})-(\ref{IV7}), and using standard 
Sobolev embeddings
 yields the result.\\
$\bullet$ {\bf Proof of i.bis} When $-t_0<s\leq t_0+1$ 
one can replace the r.h.s. of (\ref{adapt})
by 
$\Vert \sigma^1\Vert_{H^{t_0+m_1\wedge n+1}_n}
\Vert \sigma^2\Vert_{H^{s+m_1\wedge n}}$ and modify subsequently (\ref{IV4}).
Similarly, one can modify (\ref{IV5}) remarking that 
$n_{0,s}((\sigma^1\sharp_n\sigma^2)_n)\lesssim 
\Vert \sigma^1\Vert_{H^{t_0+1}_n}\Vert \sigma^2\Vert_{H^{s+n}}$.
Remarking also that in (\ref{IV6}) on can  replace
$\vert u\vert_{H^{m_1+m_2+t_0-m_1\wedge n}}$
by $\vert u\vert_{H^{s+m_1+m_2-n-1}}$ and that (\ref{IV7}) can be replaced
by
$$
	\vert \Op(\tau^7_1)u\vert_{H^s}\lesssim
	N_{0,t_0+n+1}^{m_1}(\sigma^1)M_d^{m_2}(\sigma^2)
	\vert u\vert_{H^{s+m_1+m_2-n-1}}
$$
if one uses Prop. \ref{propII5}.ii (with $r'=n+1$) rather than Prop. \ref{propII5}.i,
one gets {\bf i.bis}.\\
{\bf Proof of i.ter} In the proof of {\bf i.}, the only estimate which is not valid
when $t_0+1< s\leq s_0+1$ is (\ref{IV1}). To give control of $\tau_1^1$ and
$\tau_1^2$ one now has to use the second estimate of Th. \ref{theoII1} instead
of the first one, whence the additional term in the final estimate.\\
{\bf Proof of ii.} As said above, the control of $\tau_2$ is deduced
from the control of $\tau_1$ by a simple permutation of $\sigma^1$
and $\sigma^2$, whence the result.\\
{\bf Proof of iii.} When $n=0$, one has $\tau_1^4=\tau^4_2$, so that the 
control of both term is not needed to estimate the commutator. We keep
the same control of $\tau_1^3$ as in the proof of {\bf i} and explain the
modifications that must be performed to control $\tau_1^j$, $j=1,2,5,6,7$
(of course, the components of $\tau_2$ are treated the same way).\\
{\bf Control of $\mathbf{\tau^1_1(x,D)}$ and $\mathbf{\tau_1^2(x,D)}$.}
Since $\sigma^2$ is a function, one can invoke Prop. \ref{propII5bis}.iv
(see also Remark \ref{remII5bis})
and Lemma \ref{lemmII0} to deduce from (\ref{blob}) that 
for all $0\leq s\leq t_0+1$,
\begin{equation}
	\label{bay1}
	\vert \tau_1^j(x,D)u\vert_{H^s}\lesssim
	\Vert \sigma^1\Vert_{H^{t_0+1}}\vert \sigma^2\vert_{H^{s+m_1}}
	\vert u\vert_{\infty}\qquad (j=1,2).
\end{equation}
{\bf Control of $\mathbf{\tau^5_1(x,D)}$ } Since $\sigma_1\sharp_0\sigma^2=\sigma^1\sigma^2$, and because $\sigma^1$ is $\ind$-regular at the origin, we
can use Prop. \ref{propII1} to get $\vert \tau_1^5(x,D)u\vert_{H^s}\lesssim
N_{\ind,s}(\sigma^1\sigma^2)\vert u\vert_\infty$. It follows easily that
for all $0\leq s\leq t_0+1$,
\begin{equation}
	\label{bay2}
	\vert \tau_1^5(x,D)u\vert_{H^s}\lesssim 
	\Vert \sigma^1\Vert_{H^{t_0+1}_{reg}}
	\vert \sigma^2\vert_{H^{s}}\vert u\vert_\infty.
\end{equation}
$\bullet$ {\bf Control of $\mathbf{\tau_1^6(x,D)}$.} Using successively 
Props. \ref{propII1} and \ref{propII2} as in {\bf i}
 one can also obtain, for all $s\leq t_0+1$,
\begin{equation}
	\label{bay3}
	\vert \Op(\tau^6_1)u\vert_{H^s}\lesssim
	n_{0,s}(\sigma^1)\vert \sigma^2\vert_\infty
	\vert u\vert_{H^{s+m_1-1}}.
\end{equation}
$\bullet$ {\bf Control of $\mathbf{\tau_1^7(x,D)}$.} Using successively 
Props. \ref{propII5}.ii (with $r'=1$) and \ref{propII2}
 one obtains, for all $0\leq s\leq t_0+1$,
\begin{equation}
	\label{bay4}
	\vert \Op(\tau^7_1)u\vert_{H^s}\lesssim
	N_{0,t_0+1}^{m_1}(\sigma^1)\vert \sigma^2\vert_\infty
	\vert u\vert_{H^{s+m_1-1}}.
\end{equation}

Point {\bf iii} of the theorem follows from (\ref{IV3}) and 
(\ref{bay1})-(\ref{bay4}).
\end{proof}

\begin{rem}
	Taylor proved in \cite{TaylorM0}
	that for all classical pseudo-differential operator $\sigma^1(x,D)$
	of order $m_1>0$ and $\sigma^2\in H^\infty(\R^d)$, 
	one has, for all $s\geq 0$, 
	$$
	\left\vert [\Op(\sigma^1),\sigma^2]u\right\vert_{H^s}
	\leq \underline{C}(\sigma^1)\left(\vert\sigma^2\vert_{W^{1,\infty}} 
	\vert u\vert_{H^{s+m_1-1}}
	+\vert \sigma^2\vert_{H^{s+m_1}}\vert u\vert_\infty\right).
	$$
	This is exactly the estimate of Th. \ref{theoIV1}.iii, 
	which also gives a description of the
	constant $\underline{C}(\sigma^1)$. The commutator estimate
	corresponding to {\bf i.bis} generalizes this result taking
	into account the smoothing effect of the Poisson bracket. It turns
	out that this estimate is not tame with respect to $m_1$, while
	the commutator estimate corresponding to {\bf i}, which is not
	\emph{stricto sensu} of Kato-Ponce type, is tame.
\end{rem}

We also give commutator estimates of Calderon-Coifman-Meyer type:
\begin{theo}
	\label{theoIV1bis}
	Let $m_1,m_2\in\R$, $n\in\N$ and $d/2<t_0\leq s_0$.
	Let 
	$\sigma^j(x,\xi)\in\Gamma^{m_2}_{s_0+m_j\wedge n+1}$
	 ($j=1,2$) be $n$-regular at the origin.\\
	For all $s\in\R$ such that	
	$-t_0<s+m_j\leq t_0+n+1$ ($j=1,2$) and $-t_0<s\leq t_0+1$, one has
	$$
	\big\vert \big[\Op(\sigma^1),\Op(\sigma^2)\big]u
	-\Op(\{\sigma^1,\sigma^2\}_n)u
	\big\vert_{H^s}
	\lesssim  \Vert \sigma^1\Vert_{H^{t_0+n+1}_n}
	\Vert \sigma^2\Vert_{H^{t_0+n+1}}
	\vert u\vert_{H^{s+m_1+m_2-n-1}}.
	$$
\end{theo}
\begin{proof}
The proof follows closely the proof of Th. \ref{theoIV1}, so that we
just mention the adaptations that have to be made.\\
$\bullet$ {\bf Control of $\mathbf{\tau^j_1(x,D)}$, $j=1,2,4,5$.}  
One just has to do as in the proof of
Th. \ref{theoIII1bis}.\\
$\bullet$ {\bf Control of $\mathbf{\tau^6_1(x,D)}$.} When gets easily
from Props. \ref{propII1} and \ref{propII2} that for all $s\leq t_0+1$,
one has $\vert \tau_1^6(x,D)u\vert_{H^s}\lesssim 
\Vert \sigma^1\Vert_{H^{t_0+1}}
\Vert \sigma^2\Vert_\infty \vert u\vert_{H^{s+m_1+m_2-n-1}}$.\\
$\bullet$ {\bf Control of $\mathbf{\tau^7_1(x,D)}$.} By Prop. \ref{propII5}.ii
(with $r'=n+1$) and Prop. \ref{propII2}, one obtains
 $\vert \tau_1^7(x,D)u\vert_{H^s}\lesssim 
\Vert \sigma^1\Vert_{H^{t_0+n+1}}
\Vert \sigma^2\Vert_\infty \vert u\vert_{H^{s+m_1+m_2-n-1}}$, 
for all $-t_0<s\leq t_0+1$.
\end{proof}

\begin{rem}
	When $n=0$ and $\sigma^2$ is a function, the estimate of 
	Th. \ref{theoIV1bis} is exactly the Calderon-Coifman-Meyer 
	estimate (\ref{intro3}), with
	extended range of validity. Th. \ref{theoIV1bis} is also more general
	in the sense that it allows $n>0$ and $\sigma^2$ to 
	be a pseudo-differential operator.
\end{rem}

When the symbols $\sigma^1$ and $\sigma^2$ are of the form
(\ref{A0}), one gets the following corollary:
\begin{cor}
	\label{coroIV1}
	Let $m_1,m_2\in\R$, $m:=m_1\wedge m_2$, 
	 and $d/2<t_0\leq s_0$.
	Let also $\sigma^j(x,\xi)=\Sigma^j(v^j(x),\xi)$ with $p_j\in\N$,
	$\Sigma^j\in C^\infty(\R^{p_j},{\mathcal M}^{m_j})$ and 
	$v_j\in H^{s_0+m\wedge n+1}(\R^d)^{p_j}$ ($j=1,2$). Assume
	moreover that $\Sigma^1$ and $\Sigma^2$ are $n$-regular
	at the origin.\\
	{\bf i.} For all $s\in\R$ such that 
	$\min\{-t_0,-t_0-m_1,-t_0-m_2\}\leq s\leq s_0+1$  the
	following estimate holds (writing $v:=(v^1,v^2)$)
	\begin{eqnarray*}
	\left\vert [\Op(\sigma^1),\Op(\sigma^2)]u
	-\Op(\{\sigma^1,\sigma^2\}_n)
	\right\vert_{H^s}\lesssim C(\vert v\vert_{W^{n+1,\infty}})
	\vert u\vert_{H^{s+m_1+m_2-n-1}}\\
	+
	C(\vert v\vert_{W^{n+1,\infty}})\big(\vert v^1\vert_{H^{t_0+1}}
	\vert v^2\vert_{H^{s_++m\wedge n}}+\vert v^2\vert_{H^{t_0+1}}
	\vert v^1\vert_{H^{s_++m\wedge n}}\big)
	\vert u\vert_{H^{m+t_0}};
	\end{eqnarray*}
	{\bf ii.} For all $s\in\R$ such that	
	$-t_0<s+m_j\leq t_0+n+1$ ($j=1,2$) and $-t_0<s\leq t_0+1$, one has
	$$
	\big\vert \big[\Op(\sigma^1),\Op(\sigma^2)\big]u
	-\Op(\{\sigma^1,\sigma^2\}_n)u
	\big\vert_{H^s}
	\lesssim  C(\vert v\vert_\infty)\vert \sigma^1\vert_{H^{t_0+n+1}}
	\vert \sigma^2\vert_{H^{t_0+n+1}}
	\vert u\vert_{H^{s+m_1+m_2-n-1}}.
	$$
\end{cor}
\begin{proof}
Writing $\sigma^j(x,\xi)=\left[\sigma^j(x,\xi)-\Sigma^j(0,\xi)\right]+\Sigma^j(0,\xi)$, the result follows from Lemma \ref{lemmA0}, Ths. \ref{theoIII1} and
 \ref{theoIV1} (for {\bf i}) and Ths. \ref{theoIII1bis} and \ref{theoIV1bis}
(for {\bf ii}).
\end{proof}

%%%%%%%%%%%%%%%%%%%%%%%%%%%%%%%%%%%%%%%%%%%%%%%%%%%%%%%%%%%%%%%%%%%%%%%%%%%%
%%%%%%%%%%%%%%%%%%%%%%%%%%%%%%%%%%%%%%%%%%%%%%%%%%%%%%%%%%%%%%%%%%%%%%%%%%%%
%%%%%%%%%%%%%%%%%%%%%%%%%%%%%%%%%%%%%%%%%%%%%%%%%%%%%%%%%%%%%%%%%%%%%%%%%%%%
\appendix
\section{Proof of Prop. \ref{propA2}}\label{appCM}
%%%%%%%%%%%%%%%%%%%%%%%%%%%%%%%%%%%%%%%%%%%%%%%%%%%%%%%%%%%%%%%%%%%%%%%%%%%%
%%%%%%%%%%%%%%%%%%%%%%%%%%%%%%%%%%%%%%%%%%%%%%%%%%%%%%%%%%%%%%%%%%%%%%%%%%%%
%%%%%%%%%%%%%%%%%%%%%%%%%%%%%%%%%%%%%%%%%%%%%%%%%%%%%%%%%%%%%%%%%%%%%%%%%%%%

Owing to (\ref{LP3})-(\ref{LP4}), we can write 
$\dsp (1-\psi(\xi))\sigma(x,\xi)=\sum_{q\geq -1}\sigma_q(x,\xi)\modj^m$, 
with $\sigma_q(x,\xi):=(1-\psi(\xi))\sigma(x,\xi)\varphi_q(\xi)\modj^{-m}$. 
Obviously, for all $q\geq -1$, $\sigma_q(x,\xi)=0$ if 
$\vert\xi\vert\geq 2^{q+1}$ or $\vert\xi\vert\leq 2^{q-1}$; it follows that 
the function $\dsp A_q(x,\xi):=\sum_{k\in\Z^d}\sigma_q(x,2^{q+1}(\xi-2k\pi))$ 
is $2\pi$-periodic with respect to $\xi$ and coincides with 
$\sigma_q(x,2^{q+1}\xi)$ in the box 
${\mathcal C}:=\{\xi\in\R^d, -\pi\leq \xi_j\leq\pi, j=1,\dots,d\}$. 
Therefore, we
can write $\sigma_q(x,2^{q+1}\xi)=A_q(x,\xi)\lambda(\xi)$, where 
$\lambda\in C_0^\infty(\R^d)$ is supported in $1/5\leq\vert\xi\vert\leq 6/5$ 
and $\lambda(\xi)=1$ in $1/4\leq\vert\xi\vert\leq 1$.\\
Expending $A_q(x,\xi)$ into Fourier series, one obtains 
$$
	A_q(x,\xi)
	=\sum_{k\in\Z^d}\frac{1}{(1+\vert k\vert^2)^{1+[d/2]}}c_{k,q}(x)
	e^{i\xi\cdot k},
$$
with
$$
	c_{k,q}(x)=(1+\vert k\vert^2)^{1+[d/2]}(2\pi)^{-d}\int_{{\mathcal C}}
	e^{-i\xi\cdot k}\sigma_q(x,2^{q+1}\xi)d\xi,
$$
so that
\begin{eqnarray*}
	\sigma_q(x,\xi)&=&A_q(x,2^{-q-1}\xi)\lambda(2^{-q-1}\xi)\\
	&=&\sum_{k\in\Z^d}\frac{1}{(1+\vert k\vert^2)^{1+[d/2]}}
	c_{k,q}(x)\lambda_k(2^{-q}\xi),
\end{eqnarray*}
where $\lambda_k(\xi):=e^{i\xi\cdot k/2}\lambda(\xi/2)$ and satisfies therefore
the properties announced in the statement of the proposition.\\
The last step is therefore to obtain the desired estimates on the Fourier
coefficients $c_{k,q}$. By integration by parts, one obtains first
$$
	c_{k,q}(x)=
	(2\pi)^{-d}\int_{{\mathcal C}}
	e^{-i\xi\cdot k}\left[(1-2^{2(q+1)}\Delta_\xi)^{1+[d/2]}\sigma_q\right]
	(x,2^{q+1}\xi)d\xi,
$$
which we can rewrite as
\begin{equation}
	\label{enfin1}
	c_{k,q}(x)=(2\pi)^{-d}
	\int_{{\mathcal C}} e^{-i\xi\cdot k}
	\sum_{\vert\alpha\vert\leq 2+2[d/2]}*_\alpha 2^{(q+1)\vert\alpha\vert}
	(\partial_\xi^\alpha\sigma_q)(x,2^{q+1}\xi)d\xi,
\end{equation}
where, here and in the following, $*_{\alpha}$ denotes some
numerical coefficient depending on $\alpha$ and whose
precise value is not important.\\
Recalling that for $q\geq 0$ (we omit the case $q=-1$ which does not raise
any difficulty), $\sigma_q(x,\xi)=\widetilde{\sigma}(x,\xi)\varphi(2^{-q}\xi)$
with $\widetilde{\sigma}(x,\xi):=(1-\psi(\xi))\sigma(x,\xi)\modj^{-m}$,
one obtains, for all $\alpha\in\N^d$,
$$
	\partial_\xi^\alpha\sigma_q(x,\xi)
	=\sum_{\alpha'+\alpha''=\alpha}*_{\alpha',\alpha''}
	\partial_\xi^{\alpha'}\widetilde{\sigma}(x,\xi)
	2^{-q\vert \alpha''\vert}(\partial_\xi^{\alpha''}\varphi)(2^{-q}\xi);
$$
it follows that
\begin{eqnarray}
	\label{enfin2}
	\lefteqn{2^{(q+1)\vert\alpha\vert}
	\partial_\xi^\alpha\sigma_q(x,2^{q+1}\xi)
	=}\\
	\nonumber
	& &\sum_{\alpha'+\alpha''=\alpha}*_{\alpha',\alpha''}
	\frac{2^{(q+1)\vert\alpha'\vert}}
	{\langle2^{q+1}\xi\rangle^{\vert\alpha'\vert}}
	\langle 2^{q+1}\xi\rangle^{\vert\alpha'\vert}
	\partial_\xi^{\alpha'}\widetilde{\sigma}(x,2^{q+1}\xi)
	(\partial_\xi^{\alpha''}\varphi)(2\xi).
\end{eqnarray}
Since $\dsp \frac{2^{(q+1)\vert\alpha'\vert}}
{\langle2^{q+1}\xi\rangle^{\vert\alpha'\vert}}\leq \cst$ on the support
of $\partial_\xi^{\alpha''}\varphi(2\cdot)$, it follows from 
(\ref{enfin1}) and 
(\ref{enfin2}) that
$$
	\vert c_{k,q}\vert_{H^s}
	\leq (2\pi)^{-d}\sum_{\vert\alpha\vert\leq2+2[d/2]}
	\sum_{\alpha'+\alpha''=\alpha}*_{\alpha',\alpha''}
	\int_{{\mathcal C}}
	\langle 2^{q+1}\xi\rangle^{\vert\alpha'\vert}
	\left\vert \partial_\xi^{\alpha'}
	\widetilde{\sigma}(\cdot,2^{q+1}\xi)\right\vert_{H^s}
	d\xi,
$$
from which the estimate on $c_{k,q}$ follows. The estimate on 
$\varphi_p(D)c_{k,q}$ is proved in a similar way.

%%%%%%%%%%%%%%%%%%%%%%%%%%%%%%%%%%%%%%%%%%%%%%%%%%%%%%%%%%%%%%%%%%%%%%%%
%%%%%%%%%%%%%%%%%%%%%%%%%%%%%%%%%%%%%%%%%%%%%%%%%%%%%%%%%%%%%%%%%%%%%%%%
\section{Proof of Lemma \ref{lemmfond}}\label{applemm}
%%%%%%%%%%%%%%%%%%%%%%%%%%%%%%%%%%%%%%%%%%%%%%%%%%%%%%%%%%%%%%%%%%%%%%%%
%%%%%%%%%%%%%%%%%%%%%%%%%%%%%%%%%%%%%%%%%%%%%%%%%%%%%%%%%%%%%%%%%%%%%%%%

Throughout this proof, we write, $\psi_p(\cdot):=\psi(2^{-p}\cdot)$, 
for all $p\in\Z$. By definition of $\stI^1$ and $\stI^2$, one has
\begin{eqnarray*}
	\lefteqn{\partial_\xi^\alpha\stI^1(\cdot,\xi)
	\partial_x^\alpha \stI^2(\cdot,\xi)=
	\sum_{q\geq -1}\partial_\xi^\alpha\left(
	\psi_{q-\tN}(D_x)\sigma^1(\cdot,\xi)\varphi_q(\xi)
	(1-\psi(\xi))\right)}\\
	& &\times \sum_{p\geq -1}\psi_{p-\tN}(D_x)
	\partial_x^\alpha\sigma^2(\cdot,\xi)\varphi_p(\xi)
	(1-\psi(\xi))\\
	&= &\sum_{p\geq -1} \sum_{h=0,\pm 1}
	\partial_\xi^\alpha\left(
	\psi_{p+h-\tN}(D_x)\sigma^1(\cdot,\xi)\varphi_{p+h}(\xi)
	(1-\psi(\xi))\right)\\
	& &\times \psi_{p-\tN}(D_x)
	\partial_x^\alpha\sigma^2(\cdot,\xi)\varphi_p(\xi)
	(1-\psi(\xi)),
\end{eqnarray*}
the last equality being a consequence of the fact that for all $p\geq -1$,
one has $\varphi_p\varphi_q=0$ for all $q\neq p,p\pm 1$.\\
Remarking that for all $p\geq -1$, the $p$-th term of the above summation
has a spectrum included in the ball $\{\vert\eta\vert\leq 2^{p+2-\tN}\}$,
and recalling that $\tN=N+3$, one deduces that
\begin{equation}
	\label{lemmfond1}
	\partial_\xi^\alpha\stI^1(\cdot,\xi)
	\partial_x^\alpha \stI^2(\cdot,\xi)=\sum_{p\geq -1} 
	\psi_{p-N}(D_x)\theta_p(\cdot,\xi)	
	\varphi_p(\xi)(1-\psi(\xi))
\end{equation}
with 
\begin{eqnarray*}
	\theta_p(\cdot,\xi)&:=& \psi_{p-N-3}(D_x)
	\partial_x^\alpha\sigma^2(\cdot,\xi)\\
	& &\times \sum_{h=0,\pm 1}
	\partial_\xi^\alpha\left(
	\psi_{p+h-N-3}(D_x)\sigma^1(\cdot,\xi)\varphi_{p+h}(\xi)
	(1-\psi(\xi))\right).
\end{eqnarray*}
We now turn to study the term 
$(\partial_\xi^\alpha \sigma^1\partial_x^\alpha\sigma^2)_I$. 
By definition, one has
\begin{eqnarray*}
	(\partial_\xi^\alpha \sigma^1\partial_x^\alpha\sigma^2)_I
	(\cdot,\xi)
	&=&
	\sum_{p\geq -1} \psi_{p-N}(D_x)
	\left(\partial_\xi^\alpha\sigma^1\partial_x^\alpha\sigma^2\right)
	(\cdot,\xi)
	\varphi_p(\xi)(1-\psi(\xi))\\
	&=&\sum_{p\geq -1} \psi_{p-N}(D_x)
	\Theta_p(\cdot,\xi)
	\varphi_p(\xi)(1-\psi(\xi)),
\end{eqnarray*}
with 
$$
	\Theta_p(\cdot,\xi):=
	\sum_{h=0,\pm 1}\partial_\xi^\alpha\left(\sigma^1(\cdot,\xi)
	\varphi_{p+h}(\xi)\right)
	\partial_x^\alpha\sigma^2(\cdot,\xi),
$$
where we used (\ref{LP3}) and the fact that
$\varphi_p\varphi_q=0$ if $\vert p-q\vert\geq 2$.\\
Decomposing $\sigma^1$ into 
$\sigma^1=\psi_{p+h-N-3}(D_x)\sigma^1+(1-\psi_{p+h-N-3}(D_x))\sigma^1$, 
one obtains 
\begin{eqnarray*}
	\Theta_p(\cdot,\xi)&=&
	\sum_{h=0,\pm 1}\partial_\xi^\alpha
	\left(\psi_{p+h-N-3}(D_x)\sigma^1(\cdot,\xi)\varphi_{p+h}(\xi)\right)
	\partial_x^\alpha\sigma^2(\cdot,\xi)\\
	&+&
	\sum_{h=0,\pm 1}\partial_\xi^\alpha\left((1-\psi_{p+h-N-3}(D_x))
	\sigma^1(\cdot,\xi)\varphi_{p+h}(\xi)\right)
	\partial_x^\alpha\sigma^2(\cdot,\xi).
\end{eqnarray*}
Remark that in the first term of the r.h.s. of the above identity, one
can replace $\partial_x^\alpha\sigma^2$ by 
$\psi_{p-N+1}(D_x)\partial_x^\alpha\sigma^2$, so that one finally gets
\begin{eqnarray}
	\nonumber
	\Theta_p(\cdot,\xi)&=&\theta_p(\cdot,\xi)\\
	\nonumber
	&+&\sum_{h=0,\pm 1}\partial_\xi^\alpha
	\left(\psi_{p+h-N-3}(D_x)\sigma^1(\cdot,\xi)\varphi_{p+h}(\xi)\right)\\
	\nonumber
	& &\times
	(\psi_{p-N+1}(D_x)-\psi_{p-N-3}(D_x))
	\partial_x^\alpha\sigma^2(\cdot,\xi)\\
	\nonumber
	&+&
	\sum_{h=0,\pm 1}\partial_\xi^\alpha((1-\psi_{p+h-N-3}(D_x))
	\sigma^1(\cdot,\xi)\varphi_{p+h}(\xi))
	\partial_x^\alpha\sigma^2(\cdot,\xi)\\
	\label{lemmfond2}
	&:=&\theta_p(\cdot,\xi)+\Theta^1_p(\cdot,\xi)+\Theta^2_p(\cdot,\xi).
\end{eqnarray}
It follows therefore from (\ref{lemmfond1}) and (\ref{lemmfond2}) that
$$
	\Op\big(
	(\partial_\xi^\alpha\sigma^1\partial_x^\alpha \sigma^2)_I-
	\partial_\xi^\alpha\stI^1\partial_x^\alpha \stI^2\big)
	=\Op\big(\Theta^1(x,\xi))+\Op(\Theta^2(x,\xi)\big),
$$
where $\dsp \Theta^j(\cdot,\xi)=\sum_{p\geq -1}\psi_{p-N}(D_x)
\Theta^j_p(\cdot,\xi)\varphi_p(\xi)(1-\psi(\xi))$, $j=1,2$.\\
Quite obviously, the symbols $\Theta^1(x,\xi)$ and $\Theta^2(x,\xi)$
satisfy the assumptions of Lemma \ref{lemmII1}. The result follows 
therefore from this lemma and the estimates
\begin{equation}
	\label{lemmfond3}
	M_d^{m_1+m_2-n-1}(\Theta^1)
	\leq \cst M^{m_1}_{n+d}(\sigma^1)
	M^{m_2}_{d,n+1}(\sigma^2)
\end{equation}
and
$$
	M_d^{m_1+m_2-n-1}(\Theta^2)
	\leq \cst M^{m_1}_{n+d,n+1}(\sigma^1)
	M^{m_2}_{d,n}(\sigma^2).
$$
We only prove the first of these two estimates, the second one being obtained
in a similar way. One easily obtains that
\begin{eqnarray*}
	\vert \Theta^1(\cdot ,\xi)\vert_\infty&\leq& \cst
	\sup_{p\geq -1}
	\vert \Theta^1_p(\cdot,\xi)\varphi_p(\xi)\vert_\infty\\
	&\leq&\cst
	\sup_{p\geq -1}
	\vert \partial_\xi^\alpha(\sigma^1(\cdot,\xi)\varphi_{p}(\xi))
	\vert_\infty\\
	& &\times
	\sup_{p\geq -1} \vert (\psi_{p-N-1}(D_x)-\psi_{p-N-3}(D_x))
	\partial_x^\alpha\sigma^2(\cdot,\xi)\varphi_p(\xi)\vert_\infty.
\end{eqnarray*}
Since for all $r\in\N$, $p\geq -1$, and $f\in\sch$, one has 
$\vert (\psi_{p+2}(D)-\psi_p(D))f\vert_\infty\leq 
\cst 2^{-pr}\vert f\vert_{W^{r,\infty}}$, it follows that 
(taking $r=n+1-\vert\alpha\vert$),
$$
	\vert \Theta^1(\cdot ,\xi)\vert_\infty
	\leq \cst \modj^{m_1+m_2-n-1} M^{m_1}_n(\sigma^1)
	M^{m_2}_{0,n}(\nabla_x\sigma^2).
$$
The derivatives of $\Theta^1$ with respect to $\xi$ can be handled in the
same way, thus proving (\ref{lemmfond3}).

\bigbreak

\noindent
{\bf Acknowledgments.} The author warmly thanks T. Alazard and G. M\'etivier
for fruitful discussions about this work, and B. Texier for his remarks
on a previous version of this work. This work was partially
supported by the ACI Jeunes Chercheuses et Jeunes Chercheuses ``Dispersion 
et nonlin\'earit\'es''.

%%%%%%%%%%%%%%%%%%%%%%%%%%%%%%%%%%%%%%%%%%%%%%%%%%%%%%%%%%%%%%%%%%%%%%%%%
%%%%%%%%%%%%%%%%%%%%%%%%%%%%%%%%%%%%%%%%%%%%%%%%%%%%%%%%%%%%%%%%%%%%%%%%%
%%%%%%%%%%%%%%%%%%%%%%%%%%%%%%%%%%%%%%%%%%%%%%%%%%%%%%%%%%%%%%%%%%%%%%%%%

\end{document}